%% file: organized_main.tex
\documentclass[11pt,reqno]{amsart}
\usepackage{graphicx} 
\usepackage{todonotes}
\usepackage{comment}
\usepackage[pageanchor=false,colorlinks=true,linkcolor=blue,citecolor=blue,urlcolor=blue]{hyperref}
\usepackage{tikz}
\usepackage[
backend=biber,
style=alphabetic,giveninits,
citestyle=ieee-alphabetic,
natbib=true,
uniquelist=false,
maxnames=10,
sorting=ynt,
url=false
]{biblatex}
\DeclareFieldFormat*{title}{\mkbibemph{#1}}

\include{macro}
\usepackage[mathlines]{lineno}
\usepackage[top=1in,inner=1in,outer=1in,bottom=1in]{geometry}
\usepackage{setspace}
\usepackage{cleveref}

\title[Gradient Flows over Convex Bodies for Learning Optimal Regularizers]{Gradient Flows over Convex Bodies for\\ Learning Optimal Regularizers}
\author[Nilava Metya]{Nilava Metya}
\address{Department of Mathematics, Rutgers University - New Brunswick}
\email{nilava.metya@rutgers.edu}

\author[Oscar Leong]{Oscar Leong}
\address{Institute of Natural Sciences and School of Mathematical Sciences, Shanghai Jiao Tong University}
\email{oleong@sjtu.edu.cn}

\author[Yong Sheng Soh]{Yong Sheng Soh}
\address{Department of Mathematics and Institute of Operations Research and Analytics, National University of Singapore}
\email{matsys@nus.edu.sg}

\date{}
\subjclass[2020]{52A41, 52A20, 90C48, 49K27, 47J35, 52B05.}

\begin{document}
\maketitle
\begin{abstract}
    We consider the question of describing gradient flows over the collection of convex bodies.  This is motivated by foundational questions arising in inverse problems.  Suppose we wish to solve an ill-posed inverse problem using regularization techniques. Can we learn an optimal choice of a convex regularizer directly from prior data? Given a numerical algorithm that computes such an optimal regularizer, is there a principled manner for describing the evolution of its iterates? The objective of this paper is to lay the foundation towards answering these questions.

    To this end, we consider a criterion of an optimal regularizer based on seeking a convex body that minimizes the expected gauge evaluation with respect to a prior data distribution subject to a volume constraint. The formulation is initially nonconvex, but admits hidden convexity by passing to its polar. We characterize the Karush–Kuhn–Tucker (KKT) conditions that certify global optimality. We also analyze the $L^2$-gradient flow as well as a minimizing movement scheme with the Hausdorff penalty. Our results show that the objective gap decays exponentially along the gradient flow and we obtain dimension-dependent rates for both the objective gap and distance to the optimal body for the Hausdorff scheme. Finally, we apply the KKT conditions to identify explicit optimal convex bodies. For data where the optimal nonconvex regularizer is the $\ell_p$-quasinorm, $p \in (0,1)$, the optimal convex regularizer is the volume-normalized $\ell_1$-ball. For data described by paired basis atoms, the optimal convex regularizer is the polar of a shifted zonotope.  

\end{abstract}


\tableofcontents

\input{sec_Intro_general}
\input{sec_Notation}
\input{sec_ConvexGeometry}

\input{sec_Geodesics}
\input{sec_KKT}

\input{sec_Flows}

\input{sec_OptReg}

\input{sec_Examples}

\input{sec_Conclusion}

\printbibliography

\end{document}

%% file: macro.tex
\usepackage{latexsym,amssymb,amsmath,amsthm,graphics,graphicx,float,psfrag, epsfig, color, enumerate}
\usepackage[T1]{fontenc}

\usepackage{url,amssymb,amsmath,amsthm,amscd,paralist,bbm,wrapfig}
\usepackage{bm}
\usepackage{algorithm}
\usepackage{algorithmic} 
\usepackage{stackrel}
\usepackage{mathtools}
\usepackage{mathrsfs}
\usepackage{enumitem}
\usepackage{enumerate}
\renewcommand{\d}{\mathrm{d}}
\newcommand{\abs}[1]{\left| #1 \right|}
\newcommand{\inner}[2]{\left\langle #1, #2 \right\rangle}
\newcommand{\st}{%
  \nonscript\;
  \ifnum\currentgrouptype=16
    \;\middle|\;
  \else
    \;|\;
  \fi
  \nonscript\;}
\newcommand{\lp}{\left(}
\newcommand{\rp}{\right)}

\newcommand{\fa}{~\forall~}
\newcommand{\sett}{\coloneqq}
\newcommand{\set}[1]{\ensuremath{\left\{#1\right\}}}
\newtheorem{thm}{Theorem}[section]
\newtheorem{lem}[thm]{Lemma}
\newtheorem{cor}[thm]{Corollary}

\newtheorem{claim}[thm]{Claim}
\newtheorem{prop}[thm]{Proposition}

\newtheorem{ex}[thm]{Example}
\newtheorem{ass}[thm]{Assumption}
\newtheorem{remark}[thm]{Remark}

\newcommand{\R}{\mathbb{R}}

\newcommand{\E}{\mathbb{E}}
\newcommand{\N}{\mathbb{N}}

\newcommand{\dsphere}{\mathrm{d}\sigma}
\newcommand{\Sn}{\mathbb{S}^{n-1}}

\newcommand{\body}[1]{\ensuremath{\mathscr K_\circ^{#1}}}

\newcommand{\cE}{\mathcal{E}}
\newcommand{\cF}{\mathcal{F}}

\newcommand{\cJ}{\mathcal{J}}
\newcommand{\cB}{\mathcal{B}}

\newcommand{\cA}{\mathcal{A}}

\newcommand{\sA}{\mathscr{A}}
\newcommand{\sK}{\mathscr{K}}

\newcommand{\vol}{\mathrm{vol}}

\newcommand{\fS}{\mathfrak{S}}

\DeclareMathOperator*{\argmin}{arg\,min}
\DeclareMathOperator*{\conv}{conv}

\DeclarePairedDelimiterX{\infdivx}[2]{(}{)}{%
  #1\;\delimsize\|\;#2%
}

\DeclarePairedDelimiter{\norm}{\lVert}{\rVert}

%% file: sec_Intro_general.tex
\section{Introduction}

Many variational problems arising in analysis, data science, and engineering can be understood as optimization instances over geometric objects.  For example, in shape and spectral optimization problems, the decision variable is typically a domain whose geometry determines the solution of an underlying partial differential equation \cite{bucur2005variational}. In geometric inverse problems as well as convex set estimation, the goal is typically to recover an unknown set from indirect or noisy observations \cite{guntuboyina2012optimal, soh2021fitting, gardner2006geometric}. Related design questions arise in robust optimization where the geometry of an uncertainty set determines which perturbations a model accommodates \cite{bertsimas2018data}.  

In this work, we are interested in studying a broad family of constrained optimization problems over the collection of convex bodies.  Some conceptual questions we will be interested in include: Given a generic shape optimization instance, how do we characterize its optimal solutions?  How do we certify that a given solution is indeed locally or globally optimal?  Suppose that these optimal solutions are described as the output of an iterative algorithm.  Is there a principled manner in which we may describe the evolution of the iterates of these algorithms?  

To formally discuss these questions, we consider the following mathematical abstraction of a fairly generic class of optimization instances over convex bodies.  Let $\mathscr{K}^n$ denote the collection of nonempty, compact, convex subsets of $\R^n$. Let $\mathcal{F} : \mathscr{K}^n \rightarrow \R\cup\set{+\infty}$ be a functional that models our objective. Let $\mathscr C \subseteq \mathscr{K}^n$ denote an admissible family of solutions (convex bodies), and let $\mathcal{F}_1,\dots,\mathcal{F}_m$ be functionals that incorporate any additional geometric constraints.  In this paper, we consider the following constrained problem 
\begin{equation}
    \min_{K \in \mathscr{K}^n}\  \mathcal{F}(K)\qquad \text{subject to}\qquad K\in\mathscr{C},\ \mathcal{F}_j(K) \leq 0,\ j=1,\dots,m. \label{eq:general-problem}
\end{equation} 
The formulation \eqref{eq:general-problem} belongs to the broader family of shape optimization problems \cite{bucur2005variational}. When the family of shapes is constrained to the space of convex bodies, convex geometry \cite{Sch:13} provides a rich toolkit for addressing questions such as describing variations in $\mathcal{F}$ while preserving convex structure, with the notion of a support function playing a central role.  As a reminder, given a convex body $K$, its support function is given by \[h_K(\bm u) \sett \sup_{\bm x \in K} \langle \bm x,\bm u\rangle,\ \quad \bm u \in \Sn.\]  
A convex body $K$ is identified by its supporting hyperplanes, and hence $h_K$ specifies $K$.  Identifying convex bodies via their support functions is a powerful idea as it provides the space $\mathscr{K}^n$ with linear structure arising from Minkowski addition
\[h_{sK_0 + tK_1} = sh_{K_0} + th_{K_1},\ \quad s,t \ge 0.\] In particular, the distances \[\delta_p (K_0, K_1) \sett \|h_{K_0} - h_{K_1}\|_{L^p(\mathbb{S}^{n-1})}\sett \left(\int_{\Sn}|h_{K_0}(\bm u) - h_{K_1}(\bm u)|^p\mathrm{d}\sigma(\bm u)\right)^{1/p},\ \quad 1 \leq p \leq \infty,\] embed the geometry of convex bodies into familiar function spaces with Minkowski interpolation specifying constant-speed geodesics.  Two specific settings that will play a prominent role in this paper are when $p=2$, which places support functions in a Hilbert space, and when $p=\infty$, which recovers the usual Hausdorff distance.  These ideas are classical in convex geometry; see \cite{Sch:13}. 

Given such representations and choices of metric, an iterative descent method to solve \eqref{eq:general-problem} typically returns descent directions that balance between minimizing the objective with the amount one moves away from the current iterate.  These dynamical approaches have been studied extensively in more general metric-spaces through what is known as the {\em minimizing movement scheme} \cite{de1993new, ambrosio2005gradient}. In particular, let $\sA$ denote the feasible set in \eqref{eq:general-problem}.  Choose an initial body $K_{\tau}^0 \in \sA$, and fix a time step $\tau > 0$. Then one recursively defines \begin{align}
    K_\tau^{k+1} \in \argmin_{K \in \sA} \set{ \mathcal{F}(K) + \frac{1}{2\tau} \delta_p(K, K_\tau^k)^2}. \label{eq:min-move-scheme}
\end{align}
For $p=2$, the update \eqref{eq:min-move-scheme} is the proximal-point or backward Euler discretization of the subgradient flow associated with $\mathcal{F}$ \cite{rockafellar1976monotone, ambrosio2005gradient}. In a general metric space, the limits of these discrete evolutions as $\tau \rightarrow 0^+$ lead to the theory of minimizing movements, curves of maximal slope, and evolution variational inequalities developed by De Giorgi \cite{de1993new} and subsequently generalized by Ambrosio, Gigli, and Savar{\'e} \cite{ambrosio2005gradient}.  One of the contributions of this paper is to lay the necessary groundwork and subsequently discuss the specializations of these ideas over the space $\mathscr{K}^n$.


To this end, we point out a series of related prior work.  Variational dynamics have been studied for spectral shape optimization \cite{bucur2011shape}, eigenvalue functionals on convex bodies \cite{holle2020gradient}, and entropy functionals on planar convex curves \cite{o2024gradient}. M\'{e}rigot and Oudet \cite{merigot2014handling} develop finite-dimensional approximations and proximal methods for variational problems with convexity constraints, including optimization over support functions. Our work develops geometric optimality conditions and convergent dynamics for a class of constrained convex-body problems, and uses this framework to study the geometry of optimal regularizers. 

\paragraph{\bf Learning optimal regularizers.}  The connection between the broad task of selecting an appropriate choice of regularization function for deployment on an inverse problem with the abstract geometric optimization problem in the language of \eqref{eq:general-problem} lies with the following simple observation:  A nonnegative, positively homogeneous convex regularizer determines the types of signals it favors, and is completely specified by its unit sublevel set, which is a convex set.  

More concretely, consider an inverse problem $\bm y = A(\bm x)$ where $A : \R^n \rightarrow \R^m$ is a forward model that corrupts an underlying signal $\bm x \in \R^n$ of interest. The {\em inverse problem} is to reconstruct $\bm x$ from measurements $\bm y$.  A variational estimator $\hat{\bm x}$ of $\bm x$ is given by \begin{align}
    \hat{\bm x} \in \argmin_{\bm x \in \R^n} \left\{\mathcal{L}(A(\bm x),\bm y) + \mathcal{R}(\bm x)\right\} \label{eq:intro-inv-problem}
\end{align} where $\mathcal{L}$ is a data-fidelity term and $\mathcal{R}$ is a regularizer that promotes structural properties expected of the unknown signal. Familiar examples include the $\ell_1$ norm for sparse vectors \cite{donoho2006compressed,candes2006robust,tibshirani1996regression}, group-sparsity penalties for signals with grouped variables \cite{simon2013sparse}, and the nuclear norm for low-rank matrices \cite{recht2010guaranteed}. The atomic-norm perspective relates the structure each penalty promotes to the geometry of its unit ball \cite{chandrasekaran2012convex}. For example, the vertices (extreme points) of the $\ell_1$ unit ball are the signed coordinate vectors $\pm \bm e_i$. After normalization to unit $\ell_1$ norm, a vector with one nonzero coordinate lies at a vertex, while a vector with two nonzero coordinates lies in the relative interior of an edge. More generally, if $\bm x$ has exactly $k \geq 1$ nonzero coordinates, then $\bm x/\|\bm x\|_1$ lies in the relative interior of the $(k-1)$-dimensional face determined by its support and signs. As such, the choice of deploying the $\ell_1$-norm simply reflects a prior belief that the underlying signals of interest are sparse.  In more contemporary settings where the underlying signal structure is harder to specify in advance, a natural alternative is to learn the regularizer from data. A growing literature takes this approach using parameterized families of convex regularizers based on dictionaries, semidefinite representations, or neural networks; see, for example, \cite{tovsic2011dictionary, soh2019learning, alberti2021learning, alberti2026learning, mukherjee2020learned, goujon2023neural}. Collectively, these approaches motivate a more fundamental question: what is the best convex regularizer for a given data source?

To address this question, one can directly relate this to a shape optimization problem through the notion of a gauge or Minkowski functional. In particular, suppose we are interested in finding the best regularizer for a given data source from the family of regularizers that are convex, continuous, non-negative, and positively homogeneous. Then any such function $\mathcal
R$ is the gauge of a convex set: \[\mathcal{R}(\bm x) = \|\bm x\|_K \sett \inf \{ r \geq 0 : \bm x \in r K\},\qquad K \sett \{ \bm x \in \R^n : \mathcal{R}(\bm x) \leq 1\}.\] Learning such a function is thus equivalent to learning a convex set, providing the model with a direct representation. Directions in which $K$ extends farther from the origin receive smaller penalties, making the directional preferences of the regularizer explicit. The dictionary learning task -- a widely used approach for learning sparse representations for data -- provides a concrete example of such a regularizer that is induced by a finite collection of atoms $D \sett [\bm d_1,\dots, \bm d_K]$ and is written as \begin{align*}
    \mathcal{R}_D(\bm x) \sett \inf\{\|\bm z\|_1 : D\bm z = \bm x\} = \|\bm x\|_{K_D},\ \qquad K_D = \mathrm{conv}\left\{\pm \bm d_1,\dots, \pm \bm d_K\right\}.
\end{align*}

This brings us to the main application of our theory. Let $\mu$ be a probability  measure on $\R^n$ representing the distribution of signals of interest.  Fix $\alpha \geq 1$. A natural population criterion for selecting a convex regularizer promoting structure in $\mu$ is \begin{equation}
    \min_{K\in \mathscr K^n} \stackrel[X\sim \mu]{}{\E}{\norm X_K^\alpha}= \int_{\R^n} \|\bm x\|_K^{\alpha} \mathrm{d}\mu(\bm x)\ \quad \text{subject to}\ \quad  \vol(K) \le 1. \label{eq:intro-primal-problem}
\end{equation} This problem seeks to find the convex body whose gauge minimizes the average $\alpha$-moment of the data distribution over a volume constraint. Such a criterion is natural when one aims to use the regularizer as a penalty to be minimized in the variational framework \eqref{eq:intro-inv-problem}. The volume normalization removes a trivial scaling ambiguity, since replacing $K$ with $cK$ for $c > 0$ scales the gauge by $1/c$. The problem also admits a statistical interpretation: as shown in prior work on optimal regularization \cite{leong2025optimal}, the same objective arises from maximum likelihood estimation over Gibbs densities proportional to $\exp(-\|\cdot\|_K^{\alpha})$. Thus, the optimizer can be interpreted as both an optimal geometric regularizer and as a data-adapted probabilistic model.

\subsection{Our contributions}

We develop a variational framework for learning and optimizing convex bodies, with an emphasis on geometric optimality conditions and convergent dynamics. Our main contributions are as follows.

\paragraph{\bf Optimality conditions, global certificates, and optimal convex bodies} We establish optimality conditions for the more general constrained problem \eqref{eq:general-problem}. We characterize the dual cone of support functions through weak$-*$ closed conic combinations of Dirac functionals. This yields KKT conditions that characterize global minimizers when the functionals and admissible sets are convex along Minkowski interpolation and suitable constraint qualification holds. Specialized to \eqref{eq:intro-polar-problem}, these conditions provide geometric certificates for the global optimality of candidate regularizers.

\paragraph{\bf Continuous flows and discrete convergence} The compactness theory developed makes the constrained energy a proper, lower semicontinuous convex functional on $L^2(\mathbb{S}^{n-1}).$ We establish well-posedness of its gradient flow, an evolution variational inequality, and an energy-dissipation identity. We show that the objective gap decays at rate $O(t^{-1})$, and the trajectory converges uniformly over support functions, hence in Hausdorff distance, to a minimizer. We also establish existence, uniqueness of each update, exponential convergence, and an exponential objective-gap bound for the $L^2$ implicit Euler scheme, together with convergence of its interpolants to the continuous flow as $\tau \downarrow 0$. For the Hausdorff-penalized minimizing movement scheme, we establish fixed-step convergence to the minimizer set, and convergence to a single body whenever the minimizer is unique.

\paragraph{\bf Convex reformulation of the optimal convex regularizer problem} We then specialize our results to the optimal convex regularizer problem \eqref{eq:intro-primal-problem}. Our first contribution is to show that the original problem formulation \eqref{eq:intro-primal-problem} is nonconvex, due to the volume constraint, but that it admits a hidden convexity that can be revealed through polarity. In particular, for a convex body $K$, let $Q = K^{\circ} \sett \{\bm x \in \R^n: \langle \bm x,\bm y\rangle \leq 1, \fa \bm y \in K\}$ denote the polar body of $K$. Then the gauge of $K$ enjoys the relation: $\|\bm x\|_K = h_{Q}(\bm x),$ while the volume of $K$ satisfies $\mathrm{vol}(K) = \mathrm{vol}(Q^{\circ}) = \frac{1}{n}\int_{\mathbb{S}^{n-1}}h_Q(\bm u)^{-n} \mathrm{d}\sigma(\bm u).$ Transforming the problem \eqref{eq:intro-primal-problem} via polar bodies yields the optimization problem \begin{equation}
    \min_{Q \in \mathscr{K}_\circ^n} \int_{\R^n} h_Q(\bm x)^{\alpha} \mathrm{d}\mu(\bm x)\ \quad\ \text{s.t.}\ \quad\ \frac{1}{n}\int_{\mathbb{S}^{n-1}}h_Q(\bm u)^{-n}\mathrm{d}\sigma(\bm u) \leq 1. \label{eq:intro-polar-problem}
\end{equation} 
under the following assumption 
\begin{ass}\label{ass:finite-energy}
\begin{equation}
(i)\ \ \mu(\set{\bm x: \inner{\bm u}{\bm x}>0}) > 0 \fa \bm u\in \Sn,\qquad \qquad (ii)\ \  \int_{\R^n} \|\bm x\|_2^\alpha \d\mu(\bm x) < \infty.
\end{equation}
\end{ass}

Under the above assumption, for example, used in \cite{luo2020polar}, the objective is strictly positive on every feasible body of finite energy, so every minimizer has volume one. We show that this problem \eqref{eq:intro-polar-problem} enjoys convexity with respect to both the constraint set and the objective functional. Moreover, we show that any solution to \eqref{eq:intro-polar-problem} induces a solution to \eqref{eq:intro-primal-problem} via polarity and vice versa. Furthermore, we show that the sublevel sets of the objective control the size of the polar body, while the polar-volume constraint prevents collapse in any direction. Together these estimates give uniform upper and strictly positive lower bounds on feasible support functions on each objective sublevel set, equicontinuity, and compactness in the uniform norm. They imply existence and uniqueness of minimizers for every $\alpha\ge 1$. We finally apply our general gradient flow theory to this problem and show that, for example, the objective satisfies a quadratic growth condition, ensuring the objective gap decays exponentially along the $L^2$-gradient flow and geometrically along its implicit Euler iterates.

\paragraph{\bf Geometric characterizations of optimal convex regularizers} 
Finally, we exploit our previously developed optimality conditions to directly obtain the optimal convex regularizer for several data distributions. For continuous data where it is known the optimal nonconvex regularizer is the $\ell_p$-quasinorm with $p \in (0,1)$, we show that the unique optimal convex set is the volume-normalized $\ell_1$-ball, establishing global optimality of the $\ell_1$-regularizer for distributions with a nonconvex optimal regularizer. For distributions supported on paired basis directions, we derive explicit formulas for the optimal polar body, a shifted zonotope, and the corresponding primal unit ball. These formulas reveal how the data weights determine the geometry and asymmetry of the optimizer. In particular, equally weighted signed coordinate directions recover a scaled $\ell_1$-ball as the optimal unit ball, identifying a setting in which a familiar sparsity-promoting regularizer is globally optimal. 

 We note that, up to normalization, the polar formulation \eqref{eq:intro-polar-problem} is the power-function case of the polar Orlicz-Minkowski problem studied in \cite{luo2020polar}. They characterize optimizers for discrete measures as polytopes through geometric comparison and volume normalization. Additionally, \cite[Section~4]{leong2025optimal-uncertainty} derived a finite-dimensional convex program for the two-dimensional case and used it to compute optimal unit balls for distributions supported on a finite set. Here we use the new KKT conditions on the cone of support functions to derive optimal solutions and certify their global optimality. In a parallel line of work, \cite{traonmilin2024} has also shown that sparsity-promoting regularizers, such as the $\ell_1$-norm, are optimal for sparse recovery under criteria derived from restricted isometry bounds for linear measurements. In contrast to these prior works, we also consider the algorithmic question of evolving an initial convex body to a minimizer of the polar formulation through the theory of minimizing movements.

\subsection{Organization}

This paper is organized as follows. We establish notation and conventions in Section \ref{sec:notation}. In Section \ref{sec:preliminaries}, we first review the necessary background on convex geometry and sublinear functions. Section \ref{sec:geodesics} studies Minkowski geodesics and convexity properties of geometric functionals and constraints. Then, in Section \ref{sec:optimality}, we develop the dual-cone of support functions which leads to our optimality theory and KKT conditions. Section \ref{sec:gradient-flows-euler} describes our general convergence theory constructing and analyzing the continuous $L^2$ subgradient flow, the $L^2$ implicit Euler scheme, and the Hausdorff distance-based minimizing-movement scheme. We then specialize these gradient flow results to the optimal convex regularizer problem in Section \ref{sec:application-to-our-problem}. Finally, Section \ref{sec:optimal-convex-reg-examples} exploits our KKT theory to derive explicit examples of optimal convex regularizers, and we conclude with a discussion of future directions in Section \ref{sec:conclusion}.

%% file: sec_Notation.tex
\section{Notation and conventions}\label{sec:notation}
We adopt the following notations and conventions in this paper. In what follows, all convex bodies are subsets of $\R^n$ where $n\ge 2$. By a {\em convex body}, we mean a nonempty, compact, convex subset of $\mathbb{R}^{n}$ \cite[Page 8]{Sch:13}. Convex bodies will be denoted with uppercase letters (e.g. $K, L, Q$). Collections of sets and functions will be denoted with a scripted typeface (e.g. $\mathscr K, \mathscr F$). Functionals on a collection of sets will be denoted with a calligraphic typeface (e.g. $\mathcal{L}, \mathcal{F}$). Vectors in $\R^n$ will be denoted using bold lowercase letters (e.g. $\bm x, \bm u \in \R^n$). $\sigma$ will denote the unnormalized surface measure on $\Sn$. We shall use the following collections throughout.
\begin{align}
H_{\bm u}^+ &\sett \set{\bm x\in \R^n: \inner{\bm u}{\bm x} > 0}, \qquad \bm u\in \R^n;\label{eq:half-space}\\
\mathcal B_p^n &\sett \set{\bm x\in \R^n : \norm{\bm x}_p \le 1},\qquad p\in[1,\infty];\label{eq:ball-def}\\
\mathscr K^n &\sett \set{K\subseteq \R^n : K \text{ is a convex body}};\label{def:bodies}\\
\mathscr K^n_\circ &\sett \set{K\in \mathscr K^n : 0 \in \operatorname{int} K};\label{def:goodbodies}\\
\mathscr F^n &\sett \set{h:\R^n\to \R \st h \text{ is subadditive and positively homogeneous}}\label{def:subl-functions};\\
\mathscr F^n_\circ &\sett \set{h\in \mathscr F^n : h(\bm u) > 0 \fa \bm u \ne 0}.\label{def:positive-subl-functions}
\end{align}

We shall discuss that, via support functions, $\mathscr K^n$ can be identified with $\mathscr F^n$, and that $\mathscr K^n_\circ$ can be identified with $\mathscr F^n_\circ$ \cite[Theorem 1.7.1]{Sch:13}, so we will use them interchangeably. Additionally, for a positively homogeneous function $h:\R^n\to \R$, we shall interchangeably use $h$ with $\bar h \sett \left.h\right|_{\Sn} : \Sn\to \R$. This is because any $\bar h:\Sn\to \R$ can be extended to a positively homogeneous function $h:\R^n\to \R$ by defining $h(\bm u) = \|\bm u\|_2\bar h(\bm u/\|\bm u\|_2)$ if $\bm u\ne 0$, otherwise $h(0) = 0$. We shall endow $\mathscr F^n$ (hence $\mathscr K^n$) with the metric $\delta_p=\|\cdot \|_{L^p} \equiv \|\cdot\|_p$ for $p\in [1,\infty]$.

%% file: sec_ConvexGeometry.tex
\section{Preliminaries} \label{sec:preliminaries}


In this section we provide the necessary background in convex geometry in order to formally develop calculus over the space of convex bodies.  An excellent reference for the material presented here can be found in the monograph by Schneider \cite{Sch:13}.


\paragraph{The support function.}  Let $K \in \mathscr K^n$.  The {\em support function} of $K$ is 
\[
h_K (\bm u) \sett \sup_{\bm x \in K} \langle \bm u,\bm x\rangle, \qquad \bm u\in \R^n.
\]
Note that $h_K$ is positively homogeneous whence it is enough to specify this function on $\Sn$. $h_K$ provides sufficient information to specify $K$. This is based on the fact that every closed convex set is specified by the intersection of its supporting halfspaces \cite[Section 1.7.1]{Sch:13}. More precisely, \[K=\bigcap_{\bm u\in \Sn} \set{\bm x\in \R^n : \inner{\bm u}{\bm x} \le h_K(\bm u)}.\] In particular, the support function in a given direction $\bm u \in \Sn$ is exactly the displacement (from the origin) of the supporting hyperplane whose normal is $\bm u$.

The support function also behaves well with respect to Minkowski addition.  Recall that the 
{\em Minkowski sum} of two closed subsets $A,B \subset \R^{n}$ is given by
\[
A+B = \set{ \bm a+\bm b : \bm a \in A, \bm b \in B } \subseteq \mathbb{R}^{n}.
\]
The support of a Minkowski sum is the linear sum of the individual support functions: given $K,L \in \mathscr{K}^n$, one has $h_{K+L} = h_{K}+h_{L}$ \cite[Theorem 1.7.5]{Sch:13}.  

\paragraph{\bf The gauge and radial function.}  The gauge and radial functions are other functions that specify convex sets.  First, given $K \in \body n$, the {\em gauge function} $\| \cdot \|_K$ is defined by
\[
\|\bm x\|_K \sett \inf \{ \lambda \geq 0 : \bm x \in \lambda K \} \ \text{ for } \ \bm x \in \mathbb{R}^{n}.
\]
The support function and the gauge function are related to each other via a suitable notion of duality.  Given a subset $K \subset \mathbb{R}^{n}$, the {\em polar} $K^\circ \subset \mathbb{R}^{n}$ is defined by
\[
K^\circ \sett \{ \bm x : \langle \bm x,\bm y \rangle \leq 1 \fa \bm y \in K \}.
\]
Given $K \in \body n$, one has that $\| \cdot \|_K = h_{K^\circ}$.  The polar operation also specifies a genuine duality over $\mathscr{K}^n_{\circ}$ in the sense that if $K \in \body n$, then $K^\circ \in \body n$, and $K^{\circ \circ} = K$ \cite[Theorem 1.6.1]{Sch:13}.  

Second, the {\em radial function} of $K \in \mathscr{K}^n_\circ$ is 
\[
\rho_{K} (\bm x) \sett \sup \{ \lambda \geq 0 : \lambda \bm x \in K \}.
\]
Then, $\rho_K(\bm x) = 1/\|\bm x\|_K = 1/h_{K^\circ}(\bm x)$.

These functions provide sufficient information to specify the underlying set: given a closed convex $K$, the 1-sublevel set of the gauge of $K$ is $K$ itself. The set $K$ can be recovered from the radial by taking its reciprocal to obtain the gauge.

\paragraph{\bf Sublinear functions.}  We say that a function $f : \mathbb{R}^{n} \rightarrow \mathbb{R}$ is {\em subadditive} if $f(\bm x+\bm y) \leq f(\bm x)+f(\bm y)$ for all $\bm x,\bm y \in \mathbb{R}^{n}$.  We say that $f : \mathbb{R}^{n} \rightarrow \mathbb{R}$ is {\em positively homogeneous} if $f(\lambda \bm x) = \lambda f(\bm x)$ for all $\lambda \geq 0$ and all $\bm x \in \mathbb{R}^{n}$. We say that $f$ is {\em sublinear} if it is both subadditive and positively homogeneous. The support function of a convex body is always sublinear. It turns out that the converse is also true -- every sublinear function can be identified as the support function of a convex set.  More precisely, if $f : \mathbb{R}^{n} \rightarrow \mathbb{R}$ is sublinear, then there is a unique convex body $K \in \mathscr{K}^n$ whose support function is $f$ \cite[Theorem 1.7.1]{Sch:13}.  In this respect, the collection $\mathscr K^n$ of convex bodies  can be viewed as being in a one-to-one correspondence with the collection $\mathscr F^n$ of sublinear functions. Moreover, the collection $\mathscr K^n_\circ$ of convex bodies containing the origin in their interior is in a one-to-one correspondence with the collection $\mathscr F^n_\circ$ of sublinear functions which are strictly positive away from the origin \cite[Theorem 1.7.1]{Sch:13}. These two correspondences will be essential for our subsequent development.  We will identify $\mathscr{K}^n$ with the space of sublinear functions $\mathscr{F}^{n}$.  Very often, certain calculations are more easily understood in $\mathscr{F}^{n}$, and the one-to-one correspondence allows us to freely switch between these spaces.

\paragraph{\bf Topology over $\mathscr K^n$.}  Let $\overline{h}_K = \left.h\right|_{\Sn}: \Sn \to \R $ be the support function restricted to the sphere. The $L^p$ distances over $\mathscr{K}^n$ is defined as follows for pairs of convex bodies $K$ and $L$:
\begin{equation*}
\delta_p(K,L) \sett \norm{\overline{h}_K - \overline{h}_L}_{L^p} \text{ for } p \in [1,\infty].
\end{equation*}
Here we define 
\begin{equation*}
\norm f_{L^p} \sett \left(\int_{\Sn} \abs{f}^p \dsphere \right)^{1/p}\text{ for } p\in[1,\infty) \quad  \text{ and } \quad \norm f_{L^\infty} \sett \sup_{\bm u \in \Sn} \abs{f(\bm u)},
\end{equation*}
where $\sigma$ denotes the unnormalized measure on $\Sn$.  The case where $p=+\infty$ coincides with the Hausdorff distance \cite[Lemma 1.8.14]{Sch:13}.  

An important property regarding the $L^p$ distances is that $(\sK^n,\delta_p)$ specifies a complete metric space for all $p \in [1,\infty]$ \cite[Theorem 1.8.3 and Section 1.8 Remark 9]{Sch:13}.  In this paper, the case $p=2$ and $p=+\infty$ will be particularly significant to us. In particular, the choice of $p=2$ allows us to naturally identify support functions as a closed convex cone of the Hilbert space $L^2(\Sn)$.

%% file: sec_Geodesics.tex
\section{Geodesics over $(\sK^n,\delta_p)$} \label{sec:geodesics}

In this section we describe the geodesic properties in the space $(\sK^n,\delta_p)$. All inequalities between the support or radial functions here are pointwise on the sphere.

First we recall some concepts. Let $(X,d)$ be a complete metric space.  We say that a curve $\{x_{\theta} \}_{\theta \in [0,1]} \subseteq X$ is a (constant speed) geodesic curve if $d(x_s,x_t) = (t-s) d(x_0,x_1)$ for all $0 \leq s \leq t \leq 1$ \cite[Definition 2.4.2]{ambrosio2005gradient}.  We say that a functional $\mathcal{F} : X \rightarrow (-\infty,+\infty]$ is {\em convex} on a curve $\{x_{\theta} \}_{\theta \in [0,1]}$ if $\mathcal{F}(x_{\theta}) \leq (1-\theta) \mathcal{F}(x_0) + \theta \mathcal{F}(x_1)$ \cite[Definition 2.4.1]{ambrosio2005gradient}.  We say that a functional $\mathcal{F}$ is $\lambda$-geodesically convex if for every $x_0,x_1 \in \mathrm{dom}(\mathcal{F})$ there exists a constant speed geodesic $\{x_{\theta} \}_{\theta \in [0,1]}$ such that $\mathcal{F}$ is $\lambda$-strongly convex on it \cite[Definition 2.4.3]{ambrosio2005gradient}, namely \[\cF(x_\theta) \le (1-\theta)\cF(x_0) + \theta\cF(x_1) - \frac\lambda 2\theta(1-\theta)d(x_0,x_1)^2.\] 
Last, we say that a subset $A \subseteq X$ is geodesically convex if for all $x_0,x_1 \in A$, there is a suitable geodesic between $x_0,x_1$ lying in $A$. 

We specialize these concepts to $(\sK^n,\delta_p)$.  Let $p \in [1,\infty]$.  Suppose $K_0,K_1 \in \sK^{n}$ are convex bodies.  Define the following linear interpolation 
\[
K_{\theta} = (1-\theta) K_0 + \theta K_1\ \text{for}\ \theta \in [0,1].
\] 
For all $p \in [1,\infty]$ and all $0 \leq s \leq t \leq 1$ one has
\[
\delta_p (K_s,K_t) = \| \overline{h}_{(1-s) K_0 + s K_1} - \overline{h}_{(1-t) K_0 + t K_1} \|_p = \| (t-s)(\overline{h}_{K_0} - \overline{h}_{K_1}) \|_p = (t-s) \| \overline{h}_{K_0} - \overline{h}_{K_1} \|_p.
\]
In particular, this means that $\{K_{\theta}\}_{\theta\in[0,1]}$ identifies a geodesic in $(\sK^n,\delta_p)$, and that $(\sK^n,\delta_p)$ is a geodesic space.  For $p \in (1,\infty)$, the linear interpolation is the {\em only} geodesic.  The reason is as follows: by Minkowski's inequality we have $
\| h_{K_1} - h_{K_0}\|_{p} = \| (h_{K_1} -  h_{\tilde{K}_\theta}) + (h_{\tilde{K}_\theta} - h_{K_0}) \|_{p} \leq \| h_{K_1} -  h_{\tilde{K}_\theta} \|_{p} + \|h_{K_0} - h_{\tilde{K}_\theta}\|_{p}$. When $p \in (1,\infty)$, we have that equality holds if and only if $h_{K_1} -  h_{\tilde{K}_\theta}$ and $h_{\tilde{K}_\theta} - h_{K_0}$ are non-negatively proportional, or either of these terms are zero.  The former case recovers the linear interpolation.

\begin{prop}
    Let $V : \mathbb{R} \rightarrow \mathbb{R}$ be convex. Then $K \mapsto V(h_K(\bm u))$ is convex over $\sK^{n}$ for any $\bm u \in \mathbb{R}^{n}$. Similarly, $K \mapsto \int V(h_K (\bm u)) \d\mu(\bm u)$ is convex for any finite non-negative Borel measure $\mu$ over $\Sn$.
\end{prop}

\begin{proof}
    Let $K_0,K_1 \in \sK^{n}$, and let $K_\theta$ be the interpolant. $V$ is convex so $V(h_{K_{\theta}}) = V((1-\theta) h_{K_0} + \theta h_{K_1})) \leq (1-\theta)V(h_{K_0}) + \theta V(h_{K_1})$. Both results follow from this.
\end{proof}

\begin{prop} \label{prop:gc_v1}
Let $\alpha \geq 1$.  The functional $K \mapsto \overline{h}_K^{-\alpha}$ is geodesically convex w.r.t. $\delta_p$ over $\body{n}$.
\end{prop}

\begin{proof}
Let $K_0, K_1 \in \body{n}$.  Let $\theta \in [0,1]$, and let $K_\theta = (1-\theta) K_0 + \theta K_1$ be the linear interpolant between $K_0$ and $K_1$.  Now note that because $K_0,K_1 \in \body n$, we also have $K_\theta \in \body{n}$.  So the functions $\overline{h}_{K_0}, \overline{h}_{K_1}, \overline{h}_{K_\theta}$ take values in $(0,\infty)$.  As such, we have $\overline{h}_{K_\theta}^{-\alpha} = ((1-\theta) h_{K_0} + \theta h_{K_1} )^{-\alpha} \leq (1-\theta) h_{K_0}^{-\alpha} + \theta h_{K_1}^{-\alpha}$.  The last inequality follows from convexity of the function $t \mapsto t^{-\alpha}$ over $t \in (0,\infty)$ whenever $\alpha \geq 1$.
\end{proof}

The last result has implications regarding the (Lebesgue) volume of a convex body.  First, the volume can be expressed in terms of the radial function \cite[Eq. (1.53)]{Sch:13}
\[\vol(K) = \frac{1}{n} \int_{\Sn} \rho_{K}^{n} \d\sigma.\]



\begin{prop}
Let $\alpha > 0$.  The functional $K \mapsto \rho_K^{-\alpha}$ is geodesically convex w.r.t. $\delta_p$ over $\body{n}$.
\end{prop}

\begin{proof}
Let $K_0, K_1 \in \body{n}$.  Let $\theta \in [0,1]$, and let $K_\theta = (1-\theta) K_0 + \theta K_1$ be the linear interpolant between $K_0$ and $K_1$. First note that  $\rho_{K_{\theta}} \geq (1-\theta) \rho_{K_0} + \theta \rho_{K_1}$. Next, observe that the function $t \mapsto t^{-\alpha}$ is convex and strictly decreasing over $(0,\infty)$. We then have
\[
\rho_{K_\theta}^{-\alpha} \leq ((1-\theta) \rho_{K_0} + \theta \rho_{K_1} )^{-\alpha} \leq (1-\theta) \rho_{K_0}^{-\alpha} + \theta \rho_{K_1}^{-\alpha}.
\]
\end{proof}

%% file: sec_KKT.tex
\section{Optimality Conditions} \label{sec:optimality}


In this section we describe optimality conditions with respect to a generic optimization instance in which we minimize an objective functional over the space of convex bodies.  More formally, consider an optimization instance of the following form:
\begin{equation} \label{eq:general_convexbodyopt}
\begin{aligned}
    \min_{K \in \mathscr{K}^{n}} \quad & \cF(K) \\ 
    \mathrm{s.t.} \quad & K\in \mathscr C \\ 
    \quad & \cF_j (K) \leq 0 \quad \text{ for all } \quad 1 \leq j \leq m
\end{aligned}
\end{equation}
In Section \ref{sec:preliminaries} we noted that the collection of convex bodies is in one-to-one correspondence with the collection of sublinear functions (\cite[Theorem 1.7.1]{Sch:13}).  As such, the optimization instance \eqref{eq:general_convexbodyopt} can be expressed in terms of support functions.  Formally, we denote
\[
\mathscr{F}^n \sett \{ f : \mathbb{R}^{n} \rightarrow \mathbb{R} \st f \text{ is sublinear} \}.
\]
The functionals $\mathcal{F}$ and $\mathcal{F}_j$ may be viewed as functionals over $\mathscr{F}^{n}$ in that for all $h \in \mathscr{F}^{n}$, we define
\[
\overline{\mathcal{F}}(h) \sett \mathcal{F}(K) \text{ where } K \text{ is such that } h_{K} = h.
\]
With an abuse of notation, we shall simply write $\overline{\mathcal{F}}$ as $\mathcal{F}$.  Under this identification, \eqref{eq:general_convexbodyopt} is equivalent to the following:
\begin{equation} \label{eq:general_supportopt}
\begin{aligned}
    \min_{h \in \mathscr{F}^{n}} \quad & \mathcal{F}(h) \\ 
    \mathrm{s.t.} \quad & h \in \mathscr{C} \\
    & \mathcal{F}_j (h) \leq 0 \quad \text{ for all } \quad 1 \leq j \leq m
\end{aligned}
\end{equation}

The spaces $\mathscr{K}^{n}$ and $\mathscr{F}^{n}$ are infinite dimensional spaces for $n\ge 2$, and hence additional care is required when discussing optimality conditions.  The language of subdifferentials provide the appropriate platform to discuss these conditions.  

More formally, let $(V,\|\cdot\|)$ be a normed space, and let $V^*$ be its topological dual.  Let $\cF:V \rightarrow \R\cup\set{+\infty}$ be a closed, convex, proper function.  The subdifferential of $\cF$ at a point $u\in \operatorname{dom} \cF$ is the set of elements of the dual satisfying \cite{ABM:14}
\[
\partial \cF(u) \sett \set{ u^* \in V^*: \cF(v) - \cF(u) \geq \inner{u^*}{v-u} \text{ for all } v\in V}.
\]
To encode an admissibility constraint such as the one in \eqref{eq:general_supportopt}, let $C \subseteq V$ be a closed, convex, nonempty set. Write $\iota_C$ for the indicator of $C$, equal to zero on $C$ and $+\infty$ on $V\setminus C$. Then \cite[Definition 9.5.2, Proposition 9.5.4]{ABM:14} gives for $u\in C$ that
\[
\partial \iota_{C}(u) = N_{C} (u) \sett \{ v^* \in V^* : \langle v^* , v-u \rangle \leq 0 ~ \forall v \in C \}.
\]



We shall work with the following specialization of the Karush-Kuhn-Tucker (KKT) conditions \cite[\S 8.3, Theorem 1, pp. 217-218]{opt-vec-space}. Let $\Omega\subseteq V$ be nonempty, open and convex. Suppose $\cF_0:\Omega \to \R\cup\set{+\infty}$ is proper, convex and lower semicontinuous relative to $\Omega$, and that $\cF_i:\Omega\to\R$ are convex and continuous for $i=1,\cdots,m$. Consider the problem
\begin{equation}\label{eq:general-opt}
    \min\cF_0(v) \qquad \text{s.t.}\qquad v\in \Omega, \cF_i(v) \le 0\fa i\in[m].
\end{equation}
The subdifferential of a proper convex function $\cF$ is taken with respect to $\Omega$ as follows. For $u\in \operatorname{dom}(\cF)$, we define $\partial_\Omega \cF(u) = \set{\xi\in V^* : \cF(v)-\cF(u) \ge \inner{\xi}{v-u} \fa v\in \Omega}$. To put things into perspective, in Section \ref{sec:optimal-convex-reg-examples}, we shall take $V=(C(\Sn),\norm\cdot_\infty), C=\mathscr F^n, \Omega = \set{f\in C(\Sn) : f>0}$ with $\cF_0=E$ on $\Omega\cap C = \mathscr F_\circ^n$ and $+\infty$ on $\Omega\smallsetminus C$, where $E$ is some energy functional on $\mathscr K^n$.

\begin{prop}\label{prop:kktconditions}
In the above setting, assume Slater's condition holds, that is, $\exists \ v_0\in \Omega$ with $\cF_0(v_0) < \infty$ and $\cF_i(v_0) < 0\fa 1\le i\le m$. A feasible point $u\in\operatorname{dom}\cF_0$ is a minimizer of \eqref{eq:general-opt} if and only if there exist $\lambda_1, \cdots, \lambda_m\in \R$ such that
\begin{itemize}
\item $0\in\partial_\Omega \cF_0(u) +\sum_{i=1}^m\lambda_i\partial_\Omega \cF_i(u)$;
\item $\lambda_i\cF_i(u)=0 \fa i\in [m]$;
\item $\lambda_i\ge 0\fa i \in [m]$;
\item $\cF_i(u) \le 0\fa i\in [m]$.
\end{itemize}
The `if' direction above holds without Slater's condition.  Moreover, if the corresponding Lagrangian $\cF_0+\sum_i\lambda_i\cF_i$ is strictly convex on $\operatorname{dom}\cF_0$, then the minimizer is unique.
\end{prop}

\paragraph{\bf Choice of base space $V$.}  We proceed to apply these ideas into our context.  First, we need to make an appropriate choice of the base normed space $V$.

The first candidate is to have $V$ be the set of $L^2$-integrable functions over the unit sphere $\Sn$
\[
L^2(\Sn) : = \left.\set{ f : \Sn\to \R \st \int |f|^2 \d\sigma < +\infty }\right/ \sim
\]
where $\sim$ identifies functions which are equal $\sigma-$a.e. The choice of norm is the $L^2$-norm.  Because $L^2(\Sn)$ is a Hilbert space, the elements in the topological dual $V^*$ are naturally identified as elements in $L^2(\Sn )$.  As we discuss later, this choice of $V$ is particularly convenient for developing the theory for gradient flows precisely because $L^2(\Sn)$ is a Hilbert space.  In particular, the theory of gradient flows is most developed when the objects reside in Hilbert spaces \cite{ambrosio2005gradient}.  The downside, however, is that the description of the normal cone to $\mathscr{F}^{n}$ is more cumbersome to describe.  

This motivates our second candidate, which is the space $C(\Sn)$ of continuous functions over $\Sn$. The accompanying choice of norm is $\|\cdot\|_\infty$.  As a reminder, sublinear functions are necessarily convex and hence continuous.  We have the inclusions
\[
\mathscr{F}^{n} \subseteq C(\Sn ) \subseteq L^2 (\Sn ).
\]
The choice of $V$ as $C(\Sn )$ is particularly appealing because the dual cone of $\mathscr{F}^{n}$ admits a simple description via Dirac functionals, which we describe below. 

\paragraph{\bf Cone of sublinear functions}  First, we note that $\mathscr{F}^n$ is a closed convex cone with linearity space \[\mathscr F^n \cap (-\mathscr F^n) = \set{\bm x\mapsto \inner{\bm a}{\bm x} : \bm a\in \R^n}.\] 
The dual cone can be characterized by the simple observation that $\inner{\epsilon_{\bm x}}{f} = \bar f(\bm x)$ where $\bar f:\R^n\to\R$ is the positive homogeneous extension of $f:\Sn\to\R$ and $\epsilon_{\bm x}$ is the measure 
\[
\epsilon_{\bm x} = \begin{cases}
\|\bm x\|_2 \delta_{\bm x/\norm{\bm x}_{2}} & \bm x\ne 0\\
0&\bm x=0
\end{cases}.
\] 
This makes subadditivity of $f$ equivalent to the fact that $\inner{\epsilon_{\bm x}+\epsilon_{\bm y}-\epsilon_{\bm x+\bm y}}{f} \ge 0\fa \bm x,\bm y\in \R^n$.


%
%

\begin{lem}
The dual cone of $\mathscr{F}^n$ in $C(\mathbb S^{n-1})$ is 
\[
(\mathscr{F}^n)^*_{C} = \overline{ {\rm cone} \set{\epsilon_{\bm x}+\epsilon_{\bm y}-\epsilon_{\bm x+\bm y} : \bm x,\bm y\in \R^n}}^{w^*}
\] 
where $\overline{\set{\cdot}}^{w^*}$ denotes the weak-$*$ closure.
\end{lem}

\begin{proof} To simplify notation, we will write $\mathscr{F}\sett\mathscr{F}^n$ and $\mathcal D \sett \overline{ {\rm cone} \set{\epsilon_{\bm x}+\epsilon_{\bm y}-\epsilon_{\bm x+\bm y} : \bm x,\bm y\in \R^n}}^{w^*}$.
First, we note that the dual space of $C(\mathbb S^{n-1})$ is the set of finite signed Radon measures along with the pairing 
\[
\inner{\mu}{f}_C = \int f\mathrm{d} \mu.
\] 
It therefore follows by definition that $\mathscr{F}^* = \set{\mu \in C(\mathbb S^{n-1})^* : \inner{\mu}{f}_C \ge0 \fa f\in \mathscr{F}}$. Since every $f\in \mathscr{F}$ is sublinear, we have 
\[
\bar f\lp \bm x+\bm y \rp \le \bar f(\bm x)+\bar f(\bm y) \qquad \fa \bm x,\bm y \in \R^n
\]
where $\bar f$ is the positive homogeneous extension of $f$ from $\Sn$ to $\R^n$. Hence $\mu_{\bm x,\bm y}\sett \epsilon_{\bm x}+\epsilon_{\bm y}-\epsilon_{\bm x+\bm y}$ satisfies $\inner{\mu_{\bm x,\bm y}}{f}_C \ge 0 \fa f\in \mathscr{F}, \bm x, \bm y \in \R^n$. This implies $\mu_{\bm x,\bm y} \in \mathscr{F}^*$, from which we have ${\rm cone}\set{\mu_{\bm x,\bm y} : \bm x,\bm y\in \R^n} \subseteq \mathscr{F}^*$. Now notice that $\mathscr F^*$ is weak$-*$ closed whence $\mathcal D \subseteq \mathscr F^*$. Indeed, \[\mathscr F^* = \bigcap_{f\in \mathscr F} \set{\mu:\inner{\mu}{f}\ge 0}\]
is closed since $\mu\mapsto \inner{\mu}{f}$ is continuous in the weak$-*$ topology so each $\set{\mu:\inner{\mu}{f}\ge 0}$ is closed.

We show that these are indeed all. This, in essence, relies on a hyperplane separation type argument. Suppose $\tilde \mu \in \mathscr{F}^* \smallsetminus \mathcal D$. Since $\mathcal D$ is closed (in the weak$-*$ sense) convex cone in $C(\mathbb S^{n-1})^*$, by \cite[Theorem 2.2]{DHZ:17} $\exists\ f\in C(\Sn)$ such that $\inner{\nu}{f}_C\ge 0\fa \nu\in \mathcal D$ and $\inner{\tilde \mu}{f}_C < 0$. The latter implies that $f\notin \mathscr{F}$ since $\tilde \mu \in \mathscr{F}^*$. The former implies that $\inner{\mu_{\bm x,\bm y}}{f}_C \ge 0\fa \bm x,\bm y\in \R^n$ whence $f$ is sublinear, that is, $f\in \mathscr{F}$.  This is a contradiction. Formally this shows $\mathscr{F}^*\smallsetminus \mathcal D = \varnothing$, which completes the proof.
\end{proof}

\begin{lem}
The dual cone of $\mathscr{F}^n$ in $L^2(\mathbb S^{n-1})$ is 
\[
(\mathscr{F}^n)^*_{L^2} = \set{g \in L^2(\mathbb S^{n-1}) : g\sigma \in (\mathcal F^n)^*_{C}}.
\]
\end{lem}

\begin{proof} 
Note that $L^2(S^{n-1})$ is self dual, that is, the dual is identified with $L^2(\mathbb S^{n-1})$ via the pairing 
\[
\inner{g}{f}_{L^2} = \int fg\mathrm{d} \sigma.
\]
It then follows from the definition of a dual cone that 
\begin{align*}
g\in (\mathscr{F})^*_{L^2} &~\iff~ \inner{g}{f}_{L^2} \ge 0\fa f\in \mathscr{F} \\
&~\iff~ \int fg\mathrm{d} \sigma \ge 0\fa f\in \mathscr{F} \\
&~\iff~ \inner{g\sigma}{f}_C \ge 0 \fa f\in \mathscr{F}\\
&~\iff~ g\sigma \in ( \mathscr{F})^*_C.
\end{align*}
\end{proof}

Now we move on to finding the Fr\'echet derivative of convex functionals of a special form, namely, \[
\mathcal{G}(h)\sett\int_{\Sn} \cJ(h)(\bm u) \d\nu(\bm u).
\]
Here $\cJ:U\to C(\Sn)$ is a pointwise convex functional on a convex open set $U\subset C(\Sn)$ that is Fr\'echet differentiable at $h\in U$, while $\nu$ is a finite positive Borel measure on $\Sn$. We shall show that the derivative of $\mathcal G$ is $D\cJ (h)^*\nu$ which is the unique finite signed Borel measure satisfying \[\int_{\Sn} \phi \d(D\cJ (h)^*\nu) = \int_{\Sn}D\cJ (h)[\phi] \d\nu \fa \phi\in C(\Sn).\]

\begin{lem}\label{lem:derivative-of-integral}
$\mathcal{G}$ is convex and Fr\'echet differentiable at $h$ with subdifferential $\partial_U \mathcal{G}(h)=\set{D\cJ (h)^*\nu}$.
\end{lem}


\begin{proof}
We can write $\mathcal G = L_\nu\circ \cJ$ where $L_\nu$ is the bounded linear functional \[L_\nu(f) = \int f\d\nu.\] Positivity of $\nu$ and pointwise convexity of $\cJ$ imply
convexity of $\mathcal G$. Therefore the Fr\'echet chain rule \cite[Proposition 1]{opt-vec-space} gives $D\mathcal G(h)[g] = D(L_\nu\circ \cJ)(h)[g] = L_\nu (D\cJ(h)[g]) = \inner{D\cJ(h)^*\nu}{g}$.
\end{proof}


\begin{ex}\label{ex:h^alpha}
Work on $U=\set{h\in C(\Sn) : \min h>0}$. Let $\cJ _\alpha (h) = h^\alpha$ (with $\alpha\ge 1$) and $\mathcal{G}_{\alpha,\mu}(h) = \inner{\mu}{\cJ_\alpha (h)}_{C(\Sn)}$ where $\mu$ is a finite nonnegative Borel measure on $\Sn$. Let $h\in U, g\in C(\Sn)$. Then 
\[
[D_U \cJ _\alpha(h)g](\bm u) = \alpha h(\bm u)^{\alpha-1}g(\bm u) \qquad \text{and} \qquad D_U \mathcal{G}_{\alpha,\mu} (h) = \alpha h^{\alpha-1} \mu.
\]
In particular, if we have a measure that is supported on finitely many points $W\sett\set{\bm u_i}_{i=1}^m \subseteq \Sn$
\[
\mu_W \sett \sum_{i=1}^m a_i\delta_{\bm u_i},
\]
and where $a_i>0$, then
\begin{equation}\label{eq:finitesupportenergy}
\mathcal{G}_{\alpha,\mu_W}(h) = \sum_{i=1}^m a_ih(\bm u_i)^\alpha, \qquad D \mathcal{G}_{\alpha,\mu_W} (h) = \alpha\sum_{i=1}^m a_i h(\bm u_i)^{\alpha-1}\delta_{\bm u_i}.
\end{equation}
\end{ex}


%


%% file: sec_Flows.tex
\section{Gradient flows and implicit Euler schemes}
\label{sec:gradient-flows-euler}

The optimality conditions in the previous section allowed us to certify when a given convex body is a minimizer to our general variational problem \eqref{eq:general_convexbodyopt}. We now consider the question of how to develop an algorithmic scheme to converge to such a minimizer. To this end, we study the minimizing movement scheme \eqref{eq:min-move-scheme} with two choices of metric. With the metric $\delta_2(\cdot,\cdot)$, each update is an implicit Euler step for the $L^2$-gradient flow of the constrained energy in the support-function representation. With the metric $\delta_{\infty}(\cdot,\cdot)$, the update penalizes the Hausdorff distance between successive iterates. To begin our analysis, we first identify conditions under which the sublevel set of our objective functional is compact in an appropriate topology, which will help yield convergence of the entire scheme.

\subsection{Compactness properties and uniqueness of minimizer}

Fix $p\in[1,\infty]$. We have $\body n \subseteq \mathscr K^n \subseteq (L^p(\Sn,\sigma),\delta_p)$. Recall that $\mathscr K^n$ is a closed convex cone in $L^p(\Sn,\sigma)$. Prior to studying the gradient flow dynamics, we first aim to show that our optimization problem of interest enjoys useful compactness and coercivity properties, which is largely due structural properties of convex sets and support functions. For example, the space of support functions enjoys the following compactness property, which can be found in \cite[Theorem 3]{vitale}.

\begin{prop}\label{prop:equivalent-metrics}
Every $L^p$ bounded sequence in $\mathscr K^n$ has a uniformly convergent subsequence with the limit in $\mathscr K^n$. Moreover, every $\delta_p$, $p\in[1,\infty]$, induces the same topology on $\mathscr K^n$.
\end{prop}

Now assume $\cE:\mathscr K^n\to(-\infty,+\infty]$ is proper and lower semicontinuous with respect to $\delta_p$. We now show that if the sublevel sets of $\cE$ are bounded, then they are also necessarily compact, which then implies $\cE$ attains a finite minimum over $S = \operatorname{dom}\cE\subseteq \mathscr K^n$.

\begin{prop}\label{prop:min-attain} Suppose for some fixed $p\in[1,\infty]$, every sublevel set $S_c = \set{h\in \mathscr K^n:\cE(h)\le c}$ is $L^p-$bounded. Then every $S_c$ is uniformly compact, and $\cE$ attains a finite minimum.
\end{prop}
\begin{proof}
By Proposition \ref{prop:equivalent-metrics}, every sequence in $S_c$ has a uniformly convergent subsequence with limit in $\mathscr{K}^n$. Uniform convergence implies $\delta_p$ convergence, so lower semicontinuity places the limit in $S_c$. Thus $S_c$ is uniformly compact. Pick any $h_0 \in \mathscr K^n$ with $\cE(h_0) <\infty$. This is possible because $\cE$ is proper. Therefore $\inf\set{\cE(h) : h\in \mathscr K^n} = \inf\set{\cE(h) : h\in S_{\cE(h_0)}}$. Since $S_{\cE(h_0)}$ is uniformly compact and $\cE$ is lower semicontinuous, $\cE$ attains a finite minimum by Weierstrass’s theorem on compact topological spaces \cite[Theorem 2.8]{convexity-banach}.
\end{proof}

As a short consequence of this result, one can also show that compactness of a sublevel set of positive support functions implies that the corresponding convex bodies enjoy a uniform inner and outer radii bounds.

\begin{cor}\label{cor:unif-ball}
Assume the hypothesis from Proposition \ref{prop:min-attain} and fix some $c\in \R$. Suppose additionally that $\varnothing \ne S_c\subseteq \mathscr K^n_\circ=\set{h\in\mathscr K^n:\min_{\mathbb S^{n-1}}h>0}$. Then there are constants $0<r_c\le R_c<\infty$ such that $r_c\cB_2^n\subseteq Q\subseteq R_c\cB_2^n\fa h_Q\in S_c$.
\end{cor}
\begin{proof}
The functionals $m:h\mapsto \min h,\ M:h\mapsto \max h$ are continuous in the uniform topology. By Proposition \ref{prop:min-attain}, each $S_c$ is compact so $m,M$ respectively attain their minimum and maximum on $S_c$. Let these numbers be $r_c, R_c$ respectively. This is equivalent to $r_c \le h_Q \le R_c \fa h_Q\in S_c$, which further means $r_c\cB_2^n\subseteq Q\subseteq R_c\cB_2^n\fa h_Q\in S_c$. If $r_c=0$ then $\exists h\in \mathscr K^n_\circ, \bm u\in \Sn$ such that $h(\bm u)=0$, which is impossible by definition of $\mathscr K^n_\circ$.
\end{proof}

The preceding results give conditions under which the functional attains its minimum. We now introduce additional structure that will allow us to establish uniqueness, study convergence toward the optimal body, and can be readily applied to our optimal convex regularization problem. In particular, the main observation is that, for some geometric objectives, scaling down a feasible body lowers its objective value. A normalization constraint then determines how far the body can be scaled. Since a dilation $K \mapsto t K$ acts on support functions by $h_K \mapsto t h_K$, we model this constraint by a functional $\Psi$ of negative homogeneity. As we will see shortly, any slack in the constraint $\Psi(h) \leq 1$ can therefore be used to scale the body and lower its objective value. The following assumption formalizes these properties together with the sublevel boundedness needed for compactness. In Section \ref{sec:application-to-our-problem}, we show that the optimal regularization problem admits a formulation satisfying all four conditions. We assume the following henceforth.

\begin{ass}\label{ass:needed-assumptions}Let $S = \operatorname{dom} \cE$. We assume the following:
\begin{enumerate}[label=(\roman*)]
\item $\cE$ is convex, proper and lower semicontinuous in $(\mathscr K^n, \delta_p)$;
\item $\cE(th)< \cE(h)\fa h\in S, t\in(0,1)$ such that $th\in S$;
\item $S=\operatorname{dom}\cE = \set{h\in \mathscr K_\circ^n:\Psi(h)\le 1}\ne \varnothing$ where $\Psi:\mathscr K_\circ^n\to (0,\infty)$ is strictly convex and $\exists q>0$ such that $\Psi(th) = t^{-q}\Psi(h)\fa h\in \mathscr K_\circ^n, t>0$.
\item There is some $p\in[1,\infty]$ such that every sublevel set $S_c$ is $L^p-$bounded, that is, $\forall\ c$, $\exists R_c$ such that $\norm h_p\le R_c\fa h\in S_c$.
\end{enumerate}
\end{ass}

Hence $\inf_{\mathscr{K}^n} \cE = \inf_S\cE$. We show that this problem has a unique solution, which we call $h^* \in S$.

\begin{lem}\label{lem:unique-min}
Under Assumption \ref{ass:needed-assumptions}, $\inf_{\mathscr K^n} \cE$ has a unique minimizer $h^*$ satisfying $\Psi(h^*)=1$.
\end{lem}
\begin{proof}
From Proposition \ref{prop:min-attain} we have that $\cE$ attains a finite minimum.  The fact that it is attained follows from that $S_c$ is uniformly compact for some choice of $c$ -- for instance, we may pick $c=\cE(h_0)$ for some $h_0\in \operatorname{dom}\cE$.  Subsequently, we have $\inf_{\mathscr K^n} \cE = \inf_{S_c} \cE$.  Now, suppose $h\ne h'\in S$ are minimizers with $c_*=\cE(h) = \cE(h')$. Consider the midpoint $g = (h+h')/2$. Since $\Psi$ is strictly convex we have that $\Psi(g) < 1$. But $\cE(g) \le (\cE(h)+\cE(h'))/2 = c_*$.  By scaling $g$ to $\tilde{g} = ag$ with $a=\Psi(g)^{1/q} \in (0,1)$ we get $\Psi\left(\tilde{g}\right)=1$ with $\cE(\tilde g) < \cE(g) \le c_*$, which contradicts the minimality of $h$ and $h'$.  Hence the minimizer $h^*$ is unique.

The same scaling argument proves that $\Psi(h^*)=1$. Suppose not, and that $h$ is the minimizer with $a\sett \Psi(h)^{1/q}<1$.  Then $\tilde h=a h$ gives $\Psi(\tilde h) = 1$ with $\cE(\tilde h) = \cE(ah) < \cE(h)$, a contradiction.
\end{proof}

\subsection{The \texorpdfstring{$L^2$}{L2}-gradient flow} \label{sec:L2-grad}

Now that we have shown compactness properties of the general problem and uniqueness of the optimal solution, we turn to studying gradient flows for converging towards such a minimizer. We begin with the $L^2$-gradient flow and work over the Hilbert space $H\sett(L^2(\Sn,\sigma),\norm{\cdot}_2)$. Identify each support function with its class in $H$ and extend $\cE$ to $H$ by setting it equal to $+\infty$ outside $\mathscr K^n$, and preserve the notation $\mathcal{E}$ for this extension. The extension is proper and convex. Its finite sublevels are compact in the uniform norm by Proposition \ref{prop:min-attain}. The embedding into $H$ is continuous, so these sublevels are closed in $H$, and hence the extension is lower semicontinuous. We write $\partial\cE$ for its Hilbert-space subdifferential.

For $\tau > 0$, define the mapping \begin{equation}\label{eq:l2-discrete-flow}
    R_\tau h = \arg\min_{y\in H} \set{\cE(y) + \frac{\|h-y\|_2^2}{2\tau}}.
\end{equation}
Start at $h_0\in \operatorname{dom}(\cE)$ and take the repeated iterations $h_{k+1} = R_\tau h_k$. Note that for fixed $\tau > 0$, this mapping is akin to taking a gradient step of size $\tau$ in gradient descent for Euclidean optimization problems. We aim to show the following results, which establish convergence of our scheme: \begin{align}
&\frac{h_k-h_{k+1}}{\tau} \in \partial \cE(h_{k+1});\label{eq:l2-euler-differential-gen}\\
&\frac{\norm{h_{k+1}-v}_2^2-\norm{h_k-v}_2^2 + \norm{h_{k+1}-h_k}_2^2}{2\tau} + \cE(h_{k+1})\le\cE(v);\label{eq:l2-euler-evi-gen}\\
&\cE(h_{k+1}) +\frac{\norm{h_{k+1}-h_k}_2^2}{\tau} \le\cE(h_k);\label{eq:l2-euler-descent-gen}\\
&\sum_{k=0}^\infty\norm{h_{k+1}-h_k}_2^2 \le \tau\left(\cE(h_0)-\cE(h^*)\right).\label{eq:l2-euler-sum-gen}\\
&0\le\cE(h_N)-\cE(h^*) \le\frac{\norm{h_0-h^*}_2^2}{2N\tau}.\label{eq:l2-euler-convergence-gen}
\end{align}

\begin{prop}\label{prop:iteration-defined-euler}
For every $h\in H, \tau>0$, the update $R_\tau h\in H$ exists uniquely. The iteration $h_{k+1}=R_\tau h_k$, with $h_0\in S$, produces iterates which satisfy \eqref{eq:l2-euler-differential-gen}--\eqref{eq:l2-euler-convergence-gen} and uniformly converges to $h^*$.
\end{prop}
\begin{proof}
\cite[Corollary 2.2.2]{ambrosio2005gradient} and strict convexity gives a unique $R_\tau h\in H$. \cite[Corollary 1.4.5]{ambrosio2005gradient} provides the subdifferential of $\cE(h) + \norm{h-h_k}_{2}^2/(2\tau)$, which is $\partial \cE(h) + (h-h_k)/\tau$. By Fermat's rule, $h_{k+1}$ satisfies $(h_k-h_{k+1})/\tau \in \partial \cE(h_{k+1})$, establishing \eqref{eq:l2-euler-differential-gen}. \cite[Corollary 4.1.3, equation (4.1.3)]{ambrosio2005gradient}, with $\lambda=0$, is exactly \eqref{eq:l2-euler-evi-gen}. Taking $v=h_{k}$ in \eqref{eq:l2-euler-evi-gen} gives \eqref{eq:l2-euler-descent-gen}. Summing this gives 
\[\sum_{k=0}^\infty \norm{h_{k+1}-h_k}_2^2 = \lim_{N\to\infty}\sum_{k=0}^N \norm{h_{k+1}-h_k}_2^2 \le \lim_{N\to\infty}\tau(\cE(h_0) - \cE(h_{N+1})) \le \tau(\cE(h_0) - \cE(h^*)).\]
This proves \eqref{eq:l2-euler-sum-gen}. Taking $v=h^*$ in \eqref{eq:l2-euler-evi-gen} and summing gives \[
2\tau\sum_{k=0}^{N-1}(\cE(h_{k+1})-\cE(h^*)) \le \sum_{k=0}^{N-1}\left(\norm{h_k-h^*}_2^2 - \norm{h_{k+1} -h^*}_2^2\right) = \norm{h_0-h^*}_2^2 - \norm{h_N-h^*}_2^2 \le \norm{h_0-h^*}_2^2.\]
Invoking $\cE(h_N)\le \cdots \le \cE(h_1) \le \cE(h_0)$ into the above inequality gives $2\tau N (\cE(h_N) - \cE(h^*)) \le \norm{h_0-h^*}_2^2$, which is exactly \eqref{eq:l2-euler-convergence-gen}. All iterates $h_k$ lie in the compact set $S_{\cE(h_0)}$. By lower semicontinuity and energy convergence, every cluster point minimizes $\cE$, and uniqueness therefore implies uniform convergence, thus $L^2$ convergence, of the entire sequence to $h^*$.
\end{proof}

The continuous $L^2$ gradient flow $g:[0,\infty)\to H$ is characterized by the differential inclusion 
\begin{equation}\label{eq:inclusion}
-\dot g(t) \in \partial \cE(g(t)) \text{ for almost every } t>0.
\end{equation}
We show that with the initial condition $g(0)=h_0\in S$, there is always a unique solution with certain regularity properties. More precisely, standard results from \cite{ambrosio2005gradient} show that there is a continuous solution $g$ with
\begin{equation}\label{eq:l2-regularity}\begin{aligned}
g&\in C([0,\infty);H)\\
\left.g\right|_{(0,T)}&\in H^1((0,T);H)\qquad \fa T<\infty\\
\dot g&\in L^2((0,\infty);H)
\end{aligned}\end{equation}
where $H^1((0,T);H) = \set{f : (0,T)\to H\st f,f'\in L^2}$. The following variational inequality, energy dissipation, and convergence rates also hold:
\begin{align}
\cE(v) \ge&\ \frac12\frac{\mathrm{d}}{\mathrm{d}t}\norm{g(t)-v}_2^2 + \cE(g(t))  &\text{for a.e. }t>0.
\label{eq:EVI-convex}\\
\int_a^t\norm{\dot g(r)}_2^2\d r =&\ \cE(g(a)) - \cE(g(t)), &\text{for all } 0\le a\leq t.
\label{eq:energy-dissipation}\\
0\leq &\ \cE(g(t))-\cE(h^*) \leq \frac{\norm{h_0 -h^*}_2^2}{2t}.
\label{eq:flow-sublinear-rate}
\end{align}

\begin{prop}[Well-posedness and convergence of the $L^2$-flow]
\label{prop:L2-flow-general}
For every $h_0\in S$, there is a unique curve $g:[0,\infty)\to H$ satisfying \eqref{eq:inclusion} and regularity conditions \eqref{eq:l2-regularity} such that $g(0)=h_0$, $g(t)\in S$. For every $v\in S$, this flow satisfies the variational inequality \eqref{eq:EVI-convex} and the energy-dissipation identity \eqref{eq:energy-dissipation}. The convergence of $g(t)$ to $h^*$ is uniform on $\Sn$ and its energy gap is given by \eqref{eq:flow-sublinear-rate}.
\end{prop}
\begin{proof}
The first step is to apply \cite[Theorem  4.0.4]{ambrosio2005gradient} with the complete space $H$, $\varphi=\cE$, and $\lambda=0$.  The lower semicontinuity and coercivity requirements hold because $\cE$ is proper, lower semicontinuous and has a minimum. In addition, Assumption 4.0.1 of \cite[Theorem  4.0.4]{ambrosio2005gradient} holds along the line segment $v_t=(1-t)v_0+tv_1$ because  $\cE$ is convex while the Hilbert polarization identity give, for every $z\in H$ and $\tau>0$, \[
\cE(v_t) + \frac{\norm{v_t-z}_2^2}{2\tau} \le (1-t) \left(\cE(v_0)+\frac{\norm{v_0-z}_2^2}{2\tau}\right) +t\left(\cE(v_1)+\frac{\norm{v_1-z}_2^2}{2\tau}\right) - \frac{t(1-t)}{2\tau}\norm{v_0-v_1}_2^2.\]
Taking $v=h^*$ in 
\cite[Theorem  4.0.4]{ambrosio2005gradient} now proves existence and uniqueness of the flow and the regularity conditions. \cite[Proposition 1.4.4, Corollary 1.4.2]{ambrosio2005gradient} identify this flow with \eqref{eq:inclusion}. \cite[Equation 4.0.13]{ambrosio2005gradient} proves \eqref{eq:EVI-convex}. \cite[Remark 1.3.3, equation (1.3.15)]{ambrosio2005gradient} proves \eqref{eq:energy-dissipation} for $a>0$.  We then need to show that \eqref{eq:energy-dissipation} extends to $a=0$. Observe by \cite[Theorem 4.0.4(ii)]{ambrosio2005gradient} that $\cE(g(t))\le\cE(h_0)$ whence by lower-semicontinuity and $g(t)\to h_0$ (as $t\to 0$) we have $\lim_{t\downarrow 0}\cE(g(t))=\cE(h_0)$. By letting $a\to 0^+$ in \eqref{eq:energy-dissipation} we have that \eqref{eq:energy-dissipation} holds for $0\le a\le t$.

Take $v=h^*$ in \eqref{eq:EVI-convex} and integrate over $[0,t]$, which gives \[
\frac12 \|h_0-h^*\|_2^2 \ge \int_0^t \left(\cE(g(z)-\cE(h^*))\right)\d z.
\]
By \eqref{eq:energy-dissipation}, $\cE(g(\cdot))$ is non-increasing. Hence the above integral is at least $t\cE(g(t))$, thus proving \eqref{eq:flow-sublinear-rate}. Now $\cE(g(t)) \le \cE(h_0)$ for every $t$, so the entire trajectory of $g$ is confined to $S_{\cE(h_0)}\subseteq S$. $S_{\cE(h_0)}$ is uniformly compact. Given any sequence $(t_j)_{j\in\N}$ of positive reals with $t_j\to \infty$, there is a subsequence of $(g(t_j))_j$ that uniformly converges to some $\bar h\in S_{\cE(h_0)}$. But uniform convergence implies convergence in $L^2$. Lower continuity of $\cE$ and \eqref{eq:flow-sublinear-rate} prove that $0\le \cE(\bar h) - \cE(h^*) \le 0$. By uniqueness of the minimizer $h^*=\bar h$.
\end{proof}

The continuous flow and the discrete flow are related by interpolating the discrete flow and taking $\tau\downarrow 0$. To be precise, let $(h_k^{(\tau)})$ be the discrete steps, with time step $\tau>0$, obtained from Proposition \ref{prop:iteration-defined-euler} and let $g$ be the flow from Proposition \ref{prop:L2-flow-general}, both starting at $h_0$. Define the interpolations \[\bar h^{(\tau)}(0) = h_0, \qquad \bar h^{(\tau)}(t) = h_k^{(\tau)} ~~\text{ if } (k-1)\tau < t \le k\tau, k\ge 1\]

\begin{prop}
The following holds for every $T>0$
\begin{equation}\label{eq:l2-euler-limit}
\lim_{\tau\downarrow 0}\sup_{0\le t\le T}\|\bar h^{(\tau)}(t) - g(t)\|_2 = 0.
\end{equation}
\end{prop}
\begin{proof}
We will apply \cite[Theorem 4.2.2]{ambrosio2005gradient} on the Hilbert space $H$ with the functional $\varphi=\cE$ and $\lambda=0$. First we recall that \cite[Assumptions 4.0.1]{ambrosio2005gradient} are satisfied. Indeed, $H$ is complete, $\cE$ is proper, lower semicontinuous and lower bounded, and the penalized objective $\cE(y) + \|h_k-y\|_2^2/(2\tau) $ is $\tau^{-1}$-strongly convex along Hilbert-space segments. \cite[Theorem 4.2.2]{ambrosio2005gradient} now gives a unique minimizing-movement limit $u$, with the entire family $\bar h^{(\tau)}$ converging to $u$ uniformly on every bounded time interval with respect to the $H-$norm. By \cite[Theorem 4.0.4 (ii), (iii)]{ambrosio2005gradient},   this limit is continuous at zero with $u(0) = h_0$, locally Lipschitz for positive times, and satisfies \eqref{eq:EVI-convex}. Its uniqueness identifies $u$ with $g$.
\end{proof}

Based on the preceding results, we see that we can generically obtain a energy gap rate of $O(1/t)$ where $t$ denotes time. One can obtain faster rates when the objective functional $\cE$ satisfies certain growth conditions. A concrete example satisfied by our optimal regularization problem to be considered in the next section is quadratic growth, namely if $\exists\ \gamma>0$ such that \begin{equation}\label{eq:qg}
\cE(h)-\cE(h^*) \ge \gamma\norm{h-h^*}_2^2\qquad \fa h\in S_{\cE(h_0)}.
\end{equation} We show here that this condition allows us to achieve fast convergence in both the discrete and continuous scheme:

\begin{lem}\label{lem:qg-estimates}
Suppose \eqref{eq:qg} holds for some $\gamma>0$. Then the discrete scheme $(h_k)$ and the continuous flow $g$ satisfy \begin{align}
\norm{h_k-h^*}_2 \le (1+\gamma\tau)^{-(k-\ell)}\norm{h_\ell-h^*}_2, && \cE(h_k) - \cE(h^*) \le (1+\gamma\tau)^{-(k-\ell)} (\cE(h_\ell)-\cE(h^*)),\label{eq:rate-euler-gen}\\
\norm{g(t)-h^*}_2 \le e^{-\gamma (t-a)}\norm{g(a)-h^*}_2, && \cE(g(t)) - \cE(h^*) \le e^{-\gamma (t-a)} (\cE(g(a))-\cE(h^*)).\label{eq:rate-flow-gen}
\end{align}
for $k\ge \ell\ge 0$ and $t\ge a\ge 0$
\end{lem}
\begin{proof} Denote $\cE^*\sett \cE(h^*)$ in this proof.
Let $\xi\in \partial\cE(h)$. Convexity and \eqref{eq:qg} give \begin{align}
&\gamma \norm{h-h^*}_2^2 \le \cE(h)-\cE^* \le \inner{\xi}{h-h^*} \le \norm\xi_2\norm{h-h^*}_2\label{eq:ip-bound}\\
\implies &\norm\xi_2^2 \ge \gamma (\cE(h)-\cE^*). \label{eq:diff-bound}
\end{align}

Write $\cE_k=\cE(h_k)$. By \eqref{eq:l2-euler-differential-gen}, $\xi\sett (h_{k-1}-h_{k})/\tau \in \partial\cE(h_{k})$. Taking $h=h_k$ gives 
\begin{align*}
\inner{h_{k-1}-h_k}{h_k-h^*} &\stackrel{\eqref{eq:ip-bound}}{\ge} \tau\gamma \norm{h_k-h^*}_2^2\\
\implies \inner{h_{k-1}-h^*}{h_k-h^*} &= \inner{h_{k-1}-h_k}{h_k-h^*} + \inner{h_{k}-h^*}{h_k-h^*} \ge (1+\tau\gamma) \norm{h_k-h^*}_2^2\\
\implies \norm{h_k-h^*}_2^2 \le (1+\tau\gamma)^{-1}&\inner{h_{k-1}-h^*}{h_k-h^*} \le (1+\tau\gamma)^{-1}\norm{h_{k-1}-h^*}_2 \norm{h_{k}-h^*}_2 \\
\implies \norm{h_k-h^*}_2 \le (1+\tau\gamma)^{-1}&\norm{h_{k-1}-h^*}_2 \le \cdots \le (1+\tau\gamma)^{-(k-\ell)}\norm{h_{\ell}-h^*}_2;\\\\
\cE_{k-1} - \cE_k &\stackrel{\eqref{eq:l2-euler-descent-gen}}{\ge} \tau \cdot \frac{\norm{h_{k-1}-h_k}_2^2}{\tau^2} = \tau\norm{\xi}_2^2 \stackrel{\eqref{eq:diff-bound}}\ge \tau\gamma (\cE_k-\cE^*)\\
\implies \cE_{k-1}-\cE^* &= (\cE_{k-1} - \cE_k ) + (\cE_{k} - \cE^* ) \ge (1+\gamma\tau)(\cE_{k} - \cE^* )\\
\implies \cE_k-\cE^* \le (1+\gamma\tau)^{-1}&(\cE_{k-1}-\cE^*) \le \cdots \le (1+\gamma\tau)^{-(k-\ell)}(\cE_{\ell}-\cE^*).
\end{align*}
This establishes \eqref{eq:rate-euler-gen}.

Write $\cE_t = \cE(g(t))$. Using $\xi=-\dot h(t) \stackrel{\eqref{eq:inclusion}}{\in} \partial \cE(h(t))$ gives 
\begin{align*}
\frac{\d}{\d t}(\norm{g(t)-h^*}_2^2) &= 2\inner{\dot g(t)}{h(t)-h^*} \stackrel{\eqref{eq:ip-bound}}\le -2\gamma\norm{g(t)-h^*}_2^2; \\
\frac{\d}{\d t}(\cE^* - \cE_t) &\stackrel{\eqref{eq:energy-dissipation}}{=} \norm{\dot g(t)}_2^2 \stackrel{\eqref{eq:diff-bound}}\ge \gamma(\cE_t-\cE^*). 
\end{align*} 
Applying Gr\"onwall's inequality establishes \eqref{eq:rate-flow-gen}.
\end{proof}

As a final application, we show that even stronger rates can be established when $\cE$ is $\lambda-$strongly convex, that is, \[\cE((1-t)h+th')\le(1-t)\cE(h)+t\cE(h')-\frac\lambda2t(1-t)\norm{h-h'}_2^2 \fa t\in[0,1].\]

\begin{lem}
\label{lem:quadratic-flow-rate}
Let $\cE$ be $\lambda-$strongly convex with $\lambda>0$. If $g(t)$ is a solution of \eqref{eq:inclusion} and $\set{h_k}$ are the iterates from \eqref{eq:l2-discrete-flow} both starting with $h_0$, then 
\begin{align}
\norm{g(t)-h^*}_2 \leq e^{-\lambda t}\norm{h_0-h^*}_2.
\label{eq:flow-distance-rate}\\
\norm{h_k-h^*}_2 \leq (1+\lambda\tau)^{-k}\norm{h_0-h^*}_2.\label{eq:discrete-distance-rate}
\end{align}
Moreover, for $0\leq a\leq t$ and $0\le \ell\le k \in \mathbb Z$,
\begin{align}
\cE(g(t))-\cE(h^*) \leq e^{-2\lambda(t-a)} \lp \cE(g(a))-\cE(h^*)\rp.
\label{eq:flow-energy-rate}\\
\cE(h_k)-\cE(h^*) \leq (1+\tau\lambda)^{-2(k-\ell)} \lp \cE(h_\ell)-\cE(h^*)\rp.\label{eq:discrete-energy-rate}
\end{align}
\end{lem}
\begin{proof}
The Hilbert polarization identity in the proof of Proposition \ref{prop:L2-flow-general} and strong convexity give $(\tau^{-1}+\lambda)$-strong convexity of the proximal objective for every base point. If $f,g$ are two solutions starting at different points $h_0,g_0\in S$, \cite[Theorem 4.0.4(iv)]{ambrosio2005gradient} gives \[
\norm{f(t)-g(t)}_2\le e^{-\lambda t}\norm{f(0)-g(0)}_2\fa t\ge0.\] 
Since $0\in \partial \cE(h^*)$, one can apply this with $f\equiv h^*$ to get \eqref{eq:flow-distance-rate}. \cite[Theorem 2.4.14]{ambrosio2005gradient} proves \eqref{eq:flow-energy-rate} for $a>0$ and the continuity-at-zero argument in the derivation of the energy identity in Proposition \ref{prop:L2-flow-general} proves it for $a=0$. The discrete versions \eqref{eq:discrete-distance-rate}, \eqref{eq:discrete-energy-rate} are consequences of \cite[Theorem 2.4.14, 3.1.6(ii)]{ambrosio2005gradient} and the proof of Lemma \ref{lem:qg-estimates}.
\end{proof}

\begin{remark}
\cite[Lemma 2.4.13]{ambrosio2005gradient} gives that a $\lambda-$strongly convex function satisfies \eqref{eq:qg} with $\gamma=\lambda/2$. However, the estimates of Lemma \ref{lem:qg-estimates} are weaker than those obtained for strongly convex functions using \cite[Theorem 2.4.14, 3.1.6]{ambrosio2005gradient}. 
\end{remark}
    
\subsection{Hausdorff minimizing movement scheme}\label{sec:hausdorff-general}
We now move on to a gradient flow for minimizing $\cE$ over $B = (C(\Sn),\|\cdot\|_\infty)$, which is a complete metric space. Let $S=\operatorname{dom}(\cE)\subseteq B$ and suppose Assumption \ref{ass:needed-assumptions} holds with the ambient space $B$ by extending $\cE$ to $B$ by $+\infty$ outside $\operatorname{dom}\cE$. Then the arguments of Propositions \ref{prop:equivalent-metrics}, \ref{prop:min-attain}, Lemma \ref{lem:unique-min} show that all sublevel sets $S_c$ are compact and $\cE$ attains its minimum at the unique minimizer $h^* \in S$. Consider the updates \begin{equation}\label{eq:hausdorff-update}
    T_\tau h \in \arg\min_{y\in B} \set{\cE(y) + \frac{\|h-y\|_\infty^2}{2\tau}}.
\end{equation}
Thus, \cite[Assumptions 2.1(a)-(c)]{ambrosio2005gradient} are satisfied with the uniform topology. Moreover if $h,g\in S$ then $\|h_t-h\|_\infty = t\|g-h\|_\infty$ where $h_t=(1-t)h+tg, t\in[0,1]$. This gives \begin{equation}\label{eq:slope-hausdorff}
\cE(h_t) + \frac{\|h-h_t\|_\infty^2}{2\tau} \le (1-t)\cE(h) + t\lp\cE(g) + \frac{\|h-g\|_\infty^2}{2\tau}\rp - \frac{t(1-t)}{2\tau}\|h-g\|_\infty^2.
\end{equation}
Thus \cite[Assumption 2.4.5]{ambrosio2005gradient} is satisfied with $\lambda=0$. By \cite[Theorem 2.4.9 and Corollary 2.4.10]{ambrosio2005gradient}, the descending slope is \begin{equation}\label{eq:slope-def}
s(h)\sett |\partial \cE|_{\infty}(h) = \sup_{g\in S\smallsetminus\set h}\frac{(\cE(h)-\cE(g))_+}{\|h-g\|_\infty}.
\end{equation}
We start at a point $h_0\in S$ with the iteration $h_{k+1}=T_\tau h_k$. Since we start from $h_0$, the iteration implies that every subsequent $h_k$ has decreasing energy $\cE$, so all iterates are restricted to the compact set $S_{\cE(h_0)}$. We know the minimizer $h^*$ of $\cE$ over $B$ must then belong to $S_{\cE(h_0)}$. Akin to the $L^2$ gradient flow case, we show that the following descent properties are satisfied by the Hausdorff minimizing movement scheme: \begin{align}
\cE(h_k) \ge \cE(h_{k+1}) + \frac{\|h_k-h_{k+1}\|_\infty^2}{\tau}\label{eq:energy-update-hausdorff}\\
s(h_{k+1}) \le \frac{\|h_k-h_{k+1}\|_\infty}{\tau}\label{eq:slope-update-hausdorff}\\
\cE(h_k) - \cE(h^*) \le \frac{D^2}{k\tau},\qquad k\ge 1 \label{eq:convergence-hausdorff}
\end{align}
where $D\sett \sup\set{\norm{h-h^*}_\infty : h \in S_{\cE(h_0)}}$.

\begin{prop}\label{prop:euler-scheme-converge}
Fix $\tau>0$. Every update in \eqref{eq:hausdorff-update} has a (not necessarily unique) solution. The iteration $h_{k+1}=T_\tau h_k$ starting with $h_0\in S$ satisfies \eqref{eq:energy-update-hausdorff}, \eqref{eq:slope-update-hausdorff}, satisfies $h_k\to h^*$ uniformly with energy rate \eqref{eq:convergence-hausdorff}.
\end{prop}
\begin{proof}Write $\norm\cdot$ to mean $\norm\cdot_\infty$. Once an iterate equals $h^*$, comparison with $h^*$ in the next proximal problem forces every subsequent iterate to equal $h^*$. So assume $\cE(h_k) > \cE(h_*)$ for every $k\ge 0$. Also let $\cE_k\sett \cE(h_k)$ and $\cE^*\sett \cE(h^*)$.

Existence of the update follows from \cite[Corollary 2.2.2]{ambrosio2005gradient}. \eqref{eq:energy-update-hausdorff} follows from the proof of \cite[Theorem 3.1.6(i)]{ambrosio2005gradient}. \eqref{eq:slope-update-hausdorff} is true by \cite[Lemma 3.1.3]{ambrosio2005gradient}. 
Iteratively summing \eqref{eq:energy-update-hausdorff} gives \begin{equation}\sum_{k=0}^\infty\frac{\norm{h_k-h_{k+1}}^2}{\tau} \le \lim_{N\to\infty} \sum_{k=0}^{N-1}(\cE_k-\cE_{k+1}) = \lim_{N\to\infty} (\cE_0-\cE_N) \le \cE_0-\cE^*.\label{eq:distance-sum-hausdorff}\end{equation}
Since the RHS is bounded, we have $\norm{h_k-h_{k+1}}\to 0$. Choosing the candidate $g=h^*$ in \eqref{eq:slope-def} and using \eqref{eq:slope-update-hausdorff} gives $0\le \cE_{k+1}-\cE^*\le D s(h_{k+1}) \le D\norm{h_k-h_{k+1}}/\tau \to 0$. This proves $\cE_{k}\to \cE^*$. Since all iterates lie in $S_{\cE_0}$, lower semicontinuity, compactness and uniqueness shows $h_k\to h_*$ uniformly. 

By definition of the update, $\cE_{k+1}-\cE^* + \norm{h_k-h_{k+1}}^2/(2\tau) \le \norm{h_k-h^*}^2/(2\tau)$. So, in particular, $\norm{h_k-h_{k+1}} \le \norm{h_k-h^*}$. Let $t=\norm{h_k-h_{k+1}}/\norm{h_k-h^*} \in [0,1]$. Take $w = (1-t)h_k + th^*$. Then $\norm{h_k-w} = t \norm{h_k-h^*} = \norm{h_k-h_{k+1}}$. The definition of the update and \eqref{eq:slope-hausdorff}, \eqref{eq:slope-update-hausdorff} respectively give
\begin{align}
\cE_{k+1} \le \cE(w) &\le (1-t)\cE_k + t\cE^*\label{eq:convex-line}\\
\cE_{k+1} - \cE^* \le \norm{h_{k+1} - h^*} \cdot s(h_{k+1}) &\le \frac{\norm{h_{k+1} - h^*} \norm{h_{k+1} - h_k}}{\tau}.\label{eq:slope-error}
\end{align}
Combining these we get 
\begin{align*}
\frac{1}{\cE_{k+1} - \cE^* } - \frac{1}{\cE_{k} - \cE^* } &= \frac{\cE_{k} - \cE_{k+1}}{(\cE_{k} - \cE^*)(\cE(h_{k+1}) - \cE^*)} \\
&\stackrel{\eqref{eq:convex-line}}{\ge} \frac{t}{\cE_{k+1}-\cE^*} \\
&\stackrel{\text{def of }t, \eqref{eq:slope-error}}{\ge} \frac{\tau \norm{h_k-h_{k+1}}}{\norm{h_k-h^*}\norm{h_{k+1} - h^*} \norm{h_{k+1} - h_k}} \\
&\stackrel{\text{def of }D}{\ge} \frac{\tau}{D^2}.
\end{align*}
Summing them proves $\cE_{k} - \cE^* \le \frac{D^2}{k\tau}$.
\end{proof}

Now let's move on to the continuous gradient flow. We will obtain the flow via interpolation. Define
\begin{equation}
\label{eq:hausdorff-interpolate}
\bar h_\tau(0)=h_0,\qquad \bar h_\tau(t)=h_k \quad\text{for }(k-1)\tau<t\le k\tau.
\end{equation}
For every sequence $\tau_j\downarrow0$ and every selection of updates, the interpolants have a subsequence converging to a curve $g\in AC^2_{\mathrm{loc}}([0,\infty);B)$. Precisely, $g:[0,\infty)\to B$ belongs to $AC^2_{\mathrm{loc}}([0,\infty);B)$ if, for every $T>0$, there exists a nonnegative $g_T\in L^2(0,T)$ such that
\[\norm{g(t)-g(a)}_\infty \le \int_a^t g_T(r)\d r \qquad\fa 0\le a\le t\le T.\] Note that the metric speed $|g'|_\infty(t)$ defined below can serve as $g_T$. For an absolutely continuous curve $g$ in $B$, its metric speed is
\[|g'|_\infty(t)=\lim_{a\to t}\frac{\norm{g(a)-g(t)}_\infty}{|a-t|},\] which exists for almost every $t$ by \cite[Theorem 1.1.2]{ambrosio2005gradient}. $g$ will satisfy the energy dissipation, differential identities and energy convergence rate \begin{align}
\cE(g(t)) +\frac12\int_a^t (|g'|_\infty(r))^2\d r +\frac12\int_a^t s^2(g(r))\d r &= \cE(g(a))&\fa 0\le a\le t\label{eq:hausdoff-energy}\\
|g'|_\infty(t)= s(g(t)),\qquad \frac{\mathrm{d}}{\mathrm{d}t}\cE(g(t))& =-s(g(t))^2& \text{for almost every } t>0,\label{eq:hausdorff-differential} \\
\cE(g(t)) - \cE(h^*) \le \frac{D^2(\cE(h_0)-\cE(h^*))}{(\cE(h_0)-\cE(h^*))t+D^2} && \text{if } \cE(h_0) > \cE(h^*)\label{eq:hausdorff-convergence-cont}
\end{align}

\begin{prop}\label{prop:hausdorff-scheme-limit}
For every sequence $\tau_j\downarrow0$ and every selection of the corresponding updates, there are a subsequence $(\tau_{j_m})$ and $g\in AC^2_{\mathrm{loc}}([0,\infty);B)$ with $g(0)=h_0$, such that, for every $T<\infty$,
\begin{equation}\label{eq:hausdorff-subsequence-limit}
    \lim_{m\to\infty}\sup_{0\le t\le T}\norm{\bar h_{\tau_{j_m}}(t)-g(t)}_\infty = 0.
\end{equation}
Every such limiting curve is a curve of maximal slope for $\cE$ satisfying \eqref{eq:hausdoff-energy}, \eqref{eq:hausdorff-differential} and $g(t) \to h^*$ uniformly as $t\to \infty$ with energy convergence rate \eqref{eq:hausdorff-convergence-cont}.
\end{prop}
\begin{proof}
The complete metric space $B$ and energy $\cE$ satisfy the lower semicontinuity, coercivity and compactness assumptions of \cite[Section 2.1]{ambrosio2005gradient} because the energy is lower semicontinuous, has a lower bound and has compact sublevels, and \eqref{eq:slope-hausdorff} verifies \cite[Assumption 2.4.5]{ambrosio2005gradient} with $\lambda=0$. \cite[Proposition 2.2.3, Corollary 2.4.11]{ambrosio2005gradient} therefore give a subsequential generalized minimizing movement starting from $h_0$, with pointwise convergence in $\norm\cdot_\infty$. Using \eqref{eq:distance-sum-hausdorff} and Cauchy-Schwarz gives for $a\le t$,
\begin{equation}
\norm{\bar h_\tau(a)-\bar h_\tau(t)}_\infty \le \sum_{i=\lceil a/\tau\rceil}^{\lceil t/\tau\rceil-1}\norm{h_i-h_{i+1}}_\infty \le \sqrt{(\cE(h_0)-\cE(h^*))(t-a+\tau)}.    
\end{equation}
This proves uniform convergence on compact time intervals. By \eqref{eq:slope-def}, $s$ is lower semicontinuous and is a strong upper gradient, by \cite[Theorem 2.4.9, Corollary 2.4.10]{ambrosio2005gradient}, and thus \cite[Theorem 2.3.3]{ambrosio2005gradient} gives the maximal-slope property and  \eqref{eq:hausdoff-energy}. \cite[Remark 1.3.3]{ambrosio2005gradient} gives \eqref{eq:hausdorff-differential} and \[-\frac{\d}{\d t}\cE(g(t)) = s(g(t))^2 \qquad \text{ for almost every }t>0.\]
If $\mathcal{E}(h_0)=\mathcal{E}(h^*)$, uniqueness of optimizer gives $g(t)=h^*$ for all $t\ge0$. Otherwise $D>0$. Since the trajectory remains in $S_{\cE(h_0)}$, the formula \eqref{eq:slope-def} gives \[e(t) \sett \cE(g(t)) - \cE(h^*) \le D s(g(t)) \implies e'(t) = \frac{\d}{\d t}\cE(g(t)) =-s(g(t))^2\le -e(t)^2/D^2.\]
\eqref{eq:hausdoff-energy} makes $e$ absolutely continuous and nonincreasing. If $e$ reaches zero, it stays zero and the curve remains at the unique minimizer. For any $t$ with $e(t)>0$, integrating on $[0,t]$ gives \[\frac{1}{e(t)} - \frac{1}{e(0)}=\int_{0}^t \frac{-e'(r)}{e(r)^2} \d r\ge \frac{t}{D^2}\implies e(t) \le \frac{D^2e(0)}{te(0)+D^2}.\]
where $e(0)=\cE(h_0)-\cE(h^*)$, thus proving \eqref{eq:hausdorff-convergence-cont}. The trajectory remains in the uniformly compact set $S_{\cE(h_0)}$. Lower semicontinuity and $\cE(g(t))\to \cE(h^*)$ imply that every limit point is a minimizer. Uniqueness therefore gives $g(t)\to h^*$ uniformly.
\end{proof}

Write $e(h)=\cE(h)-\cE(h^*), e_k=e(h_k)$, and $e(t)=e(g(t))$ for the energy gaps, where $g$ is any limit curve from Proposition~\ref{prop:hausdorff-scheme-limit}. The baseline energy bounds are \eqref{eq:convergence-hausdorff} and
\eqref{eq:hausdorff-convergence-cont}.

For $L_0=s(h_0)<\infty$, the step gap bounds are
\begin{equation}\label{eq:hausdorff-finite-slope}
\norm{h_{k+1}-h_k}_\infty\le\tau L_0,\qquad \norm{\bar h_\tau(t)-\bar h_\tau(a)}_\infty \le L_0(|t-a|+\tau).
\end{equation}
Suppose, in addition, that $\exists\ c>0, p>2$ satisfying
\begin{equation}\label{eq:hausdorff-power-growth}
e(h)\ge c\norm{h-h^*}_\infty^p \qquad\text{for all }h\in S_{\cE(h_0)}.
\end{equation}
For $e_0>0$, put $\rho=1-2/p$ and $\beta=c^{2/p}$, and set
\[G_\tau(k)= \left(e_0^{-\rho}+ \frac{\rho\beta k\tau}{1+\beta\tau e_0^\rho}\right)^{-1/\rho}, \qquad G(t)=\left(e_0^{-\rho}+\rho\beta t\right)^{-1/\rho}.\]
We will show the refined bounds
\begin{align}
 e_k&\le G_\tau(k),&
 \norm{h_k-h^*}_\infty&\le c^{-1/p}G_\tau(k)^{1/p},\label{eq:hausdorff-growth-rates-discrete}\\
 e(t)&\le G(t),&
 \norm{g(t)-h^*}_\infty&\le c^{-1/p}G(t)^{1/p}.\label{eq:hausdorff-growth-rates-cont}
\end{align}

\begin{lem}\label{lem:convex--rho}
$1-(1+z)^{-\rho}\ge\rho z/(1+z)$ for $z\ge0$ and $0<\rho\le 1$.
\end{lem}
\begin{proof}
Concavity of $x\mapsto x^{1-\rho}$ gives $(1+z)^{1-\rho} - 1 \le z(1-\rho)$. Rearranging proves the inequality.
\end{proof}

\begin{prop}\label{prop:hausdorff-growth-rates}
Under Assumption~\ref{ass:needed-assumptions}, if $L_0<\infty$, then \eqref{eq:hausdorff-finite-slope} holds and every limit curve is $L_0-$Lipschitz. Under Assumption \eqref{eq:hausdorff-power-growth}, the iterates and every limit curve satisfy \eqref{eq:hausdorff-growth-rates-discrete}, \eqref{eq:hausdorff-growth-rates-cont} when $e_0>0$. If $e_0=0$, both are stationary at $h^*$.
\end{prop}

\begin{proof}
Let $b_k=\norm{h_{k+1}-h_k}_\infty$. By \eqref{eq:slope-def}, \eqref{eq:energy-update-hausdorff}, and \eqref{eq:slope-update-hausdorff},
\begin{equation}\label{eq:discrete-descent}
\frac{b_k^2}{\tau}\le e_k-e_{k+1}\le s(h_k)b_k, \qquad s(h_{k+1})\le\frac{b_k}{\tau}.
\end{equation}
If $b_k>0$, divide the first two bounds by $b_k$. If $b_k=0$, the next-step slope bound gives $s(h_{k+1})=0$. In either case, $s(h_{k+1})\le b_k/\tau\le s(h_k)$. Induction gives $s(h_k)\le L_0$ for every $k$.  This gives $b_k \le \tau L_0$, thus proving the first part of \eqref{eq:hausdorff-finite-slope}. Consequently triangle inequality gives $\norm{\bar h_\tau(t)-\bar h_\tau(a)}_\infty = \norm{h_{\lceil t/\tau\rceil}-h_{\lceil a/\tau\rceil}}_\infty \le L_0\tau\abs{\lceil t/\tau\rceil-\lceil a/\tau\rceil} \le L_0 (\abs{t-a} + \tau)$, thus proving the second inequality in \eqref{eq:hausdorff-finite-slope}. Passing to a locally uniform limit $\tau\downarrow 0$ gives the Lipschitz constant $L_0$.

For $h\in S_{\cE(h_0)}$ with $e(h)>0$, using $g=h^*$ in \eqref{eq:slope-def} and applying \eqref{eq:hausdorff-power-growth} gives \begin{equation}\label{eq:slope-inequality}
s(h) \stackrel{\eqref{eq:slope-def}}{\ge} \frac{e(h)}{\norm{h-h^*}_\infty} \stackrel{\eqref{eq:hausdorff-power-growth}}{\ge} c^{1/p}e(h)^{1-1/p} = \sqrt{\beta e(h)^{1+\rho}} .
\end{equation}
Therefore we get 
\begin{align}
e_k-e_{k+1}\stackrel{\eqref{eq:discrete-descent}}\ge b_k^2/\tau \stackrel{\eqref{eq:discrete-descent}}\ge s(h_{k+1})^2\tau \stackrel{\eqref{eq:slope-inequality}}\ge \tau c^{2/p} e_{k+1}^{2(1-1/p)} = \beta\tau e_{k+1}^{1+\rho}\label{eq:discrete-e-gap}\\
e'(t) \stackrel{\eqref{eq:hausdorff-differential}}= -s(g(t))^2 \stackrel{\eqref{eq:slope-inequality}}\le -\beta e(t)^{1+\rho}
 \quad\text{for almost every }t>0.\label{eq:cont-e-gap}
\end{align}
While the gaps are positive, we have
\[
 e_{k+1}^{-\rho}-e_k^{-\rho} = e_{k+1}^{-\rho} \left(1-\left(\frac{e_k}{e_{k+1}}\right)^{-\rho}\right) \stackrel{\eqref{eq:discrete-e-gap}}\ge \frac{1-\left(1+\beta\tau e_{k+1}^\rho\right)^{-\rho}}{e_{k+1}^\rho}\stackrel{\text{Lem }\ref{lem:convex--rho}}\ge\frac{\rho\beta\tau}{1+\beta\tau e_{k+1}^{\rho}}
 \ge \frac{\rho\beta\tau}{1+\beta\tau e_0^\rho}.
\]
Summing the above inequality from $0$ to $k-1$ proves the first inequality in \eqref{eq:hausdorff-growth-rates-discrete}. Similarly, integrating $(e^{-\rho})'\ge\rho\beta$ from \eqref{eq:cont-e-gap} gives the first inequality in \eqref{eq:hausdorff-growth-rates-cont}. The condition $\norm{h-h^*}_\infty \le e(h)^{1/p}c^{-1/p}$ gives the second inequalities in \eqref{eq:hausdorff-growth-rates-discrete}, \eqref{eq:hausdorff-growth-rates-cont}. If a gap reaches zero, monotonicity of the energy and uniqueness of the minimizer imply stationarity.
\end{proof}

%% file: sec_OptReg.tex
\section{Application: Optimal Regularization}\label{sec:application-to-our-problem}

In the prior sections, our focus was to develop the language of gradient flows over the space of convex bodies for a broad family of objective functionals.  In this section, we specialize these developments to the specific task of optimal regularization.  More specifically, recall the conceptual task of seeking an {\em optimal convex regularizer} can be expressed as the following shape constrained optimization instance \cite{leong2025optimal,leong2025optimal-uncertainty}:
\begin{equation}\label{eq:original-problem}
\min_{K\in \mathscr K^n}\stackrel[X\sim \mu]{}{\E} \norm{X}_{K}^\alpha \qquad \text{ s.t. }  \qquad \vol(K)\le 1.
\end{equation} Here $\mu$ is a Borel probability measure on $\R^n$ that models the distribution of signals of interest. We now invoke Assumption \ref{ass:finite-energy} namely  \[(i)\ \ \mu(H^+_{\bm u}) > 0 \fa \bm u\in \Sn,\qquad \qquad (ii)\ \  \int_{\R^n} \|\bm x\|_2^\alpha \d\mu(\bm x) < \infty.\]  
The first condition corresponds to the open-hemisphere condition on the weighted angular measure in Assumption \ref{ass:assumption-nu}, as in the spherical setting of \cite{luo2020polar}, and lets us reduce the search space of \eqref{eq:original-problem} down to $\body n$. The second condition is necessary and sufficient for finiteness of the objective over $\mathscr K_\circ^n$, by Proposition \ref{prop:finite-objective}.

\subsection{Reformulation through hidden convexity}

We now make two important observations about problem \eqref{eq:original-problem}. The first is that for certain feasible convex bodies, the objective value is infinite. In particular, for convex sets that do not contain the origin in its interior, the objective value is infinite. This allows us to reduce the search space from $\mathscr K^n$ to $\mathscr K^n_\circ$. The second observation is that the constraint set, even with this additional reduction from $\mathscr K^n$ to $\mathscr K^n_\circ$, is not convex under Minkowski interpolation. We resolve this challenge by showing that polarity allows us to reformulate \eqref{eq:original-problem} to an equivalent problem that is convex in both the objective and constraint set.

We formalize these two observations. For the first, we show that any feasible convex body without the origin in their interior has infinite objective value.

\begin{prop}
Let $\alpha>0$ and suppose $\mu$ satisfies Assumption \ref{ass:finite-energy}(i). Then $\E_{X\sim \mu} \norm{X}_{K}^\alpha = +\infty$ for any $K \in \mathscr K^n\smallsetminus \body n$. 
\end{prop}
\begin{proof}
By the separating hyperplane theorem, $\exists \bm v\in \Sn$ with $\inner{\bm v}{\bm x} \le 0 \fa \bm x\in K$. Then for all $t>0$ we have $\inner{\bm v}{\bm x} \le 0 \fa \bm x\in tK$. So if $\bm y\in H^+ \sett \set{\bm z\in \R^n : \inner{\bm v}{\bm z} > 0}$ then $\norm{\bm y}_{K} = +\infty$ because $tK \cap H^+ = \varnothing \fa t>0$. Since $\mu(H^+) > 0$, we have $\E_{X\sim \mu}  \norm{X}_{K}^\alpha = +\infty$. 
\end{proof}

Subsequently, the previous result allows us to reduce \eqref{eq:original-problem} to the following
\begin{equation}\label{eq:first-reform}
\min_{K\in \mathscr K^n_\circ}\stackrel[X\sim \mu]{}{\E} \norm{X}_{K}^\alpha \qquad \text{ s.t. } \qquad \vol(K)\le 1.
\end{equation} 
Even though we restricted our search space to $\body n$, the question of finiteness of the objective still remains. This is where the necessity of Assumption \ref{ass:finite-energy}$(ii)$ comes in. 
\begin{prop}\label{prop:finite-objective}
Let $\alpha>0$. Suppose $\mu$ is a Borel probability measure on $\R^n$ and $K\in \mathscr K_\circ^n$. Then $\E_{X\sim \mu}\norm{X}_{K}^\alpha < \infty$ iff $\E_{X\sim \mu}\norm{X}_{2}^\alpha < \infty$.
\end{prop}
\begin{proof}
Since $0\in\operatorname{int}K$ and $K$ is compact (the same argument in the proof of Corollary \ref{cor:unif-ball}), $\exists R\ge r>0$ such that $r\mathcal B_2^n \subseteq K \subseteq R\mathcal B_2^n$. In other words, $R^{-1}\norm{\bm x}_2 \le \|\bm x\|_K \le r^{-1}\norm{\bm x}_2\fa \bm x\in \R^n$. This implies $R^{-\alpha}\E_{X\sim \mu} \|X\|_2^\alpha  \le \E_{X\sim \mu} \|X\|_K^\alpha \le r^{-\alpha}\E_{X\sim \mu} \|X\|_2^\alpha$. Hence $\E_{X\sim \mu} \|X\|_K^\alpha$ is finite iff $\E_{X\sim \mu} \|X\|_2^\alpha$ is finite.
\end{proof}

While any feasible convex body in $\body n$ now has finite objective value, a final obstacle remains. Namely, the constraint set can be shown to be not convex under Minkowski interpolation. 
\begin{claim}\label{cl:nonconvex-constraint}
For $n\ge2$, $\set{K\in \mathscr{K}_{\circ}^n:\vol(K)\le1}$ is not convex under Minkowski interpolation.
\end{claim}
\begin{proof}
Let $I\sett [-1/2,1/2]$ and $K\sett [-1,1]\times[-1/4,1/4]\times I^{n-2},\ L \sett [-1/4,1/4]\times[-1,1]\times I^{n-2}$. Both have volume one, whereas
\[\frac12K+\frac12L=[-5/8,5/8]^2\times I^{n-2}\]
has volume $25/16>1$.
\end{proof}
This obstruction can be circumvented by parameterizing the problem via its polar, in which case the objective as well as the constraints may be viewed as being convex. Recall that for a convex body $K \in \mathscr K^n_\circ$, its polar $Q = K^{\circ} = \{\bm x : \langle \bm x , \bm y \rangle \leq 1,\ \bm y \in K\}.$ The polar body enjoys the relation $h_Q = \|\cdot\|_K = 1/\rho_K$. Exploiting this relation along with the volume formula $\mathrm{vol}(K) = \frac{1}{n}\int_{\Sn} \rho_K^{n}\mathrm{d}\sigma$ shows the following: 
\begin{prop}\label{prop:constr-convex}
    $\{Q \in \mathscr K^n_\circ : \mathrm{vol}(Q^{\circ}) \leq 1\}$ is convex under Minkowski interpolation.
\end{prop}
\begin{proof}This follows from Proposition \ref{prop:gc_v1} and noting that \[\mathrm{vol}(Q^{\circ}) = \frac{1}{n} \int_{\Sn} \rho_{Q^{\circ}}^{n} \d\sigma = \frac{1}{n} \int_{\Sn} h_{Q}^{-n} \d\sigma.\]
\end{proof}

This result shows that when imposing the constraint for the polar body, the constraint set is convex. We now reformulate the problem \eqref{eq:first-reform} using the polarity transformation $K \mapsto K^{\circ} = Q$ in both the objective and constraint set.  

\subsubsection{Polar form and radial integral representation of optimal regularization} First, recall that the polar operation is an involution on $\mathscr K_\circ^n$ \cite[Theorem 1.6.1]{Sch:13} with $h_{K^\circ} = \norm{\cdot}_{K}$, $\rho_{K} = 1/h_{K^\circ}$ and $\vol(K) = n^{-1}\int_{\Sn}h_{K^\circ}^{-n}\d\sigma$ for any $K\in \body n$. Writing $Q=K^\circ$,  \eqref{eq:first-reform} becomes
\begin{equation}
\label{eq:second-reform}
\min_{Q\in \mathscr K^n_\circ}\stackrel[X\sim \mu]{}{\E} h_Q(X)^\alpha \qquad \text{ s.t. } \qquad \frac1n\int_{\Sn}h_{Q}^{-n}\d\sigma\le 1.
\end{equation}

Projecting the data measure onto $\Sn$ with a radial weight yields a spherical integral representation of the above polar objective. This new radial measure gives a cleaner equivalent reformulation of Assumption \ref{ass:finite-energy}.

Given a Borel $A\subset \Sn$ we define
\begin{equation}\label{eq:def-nu}\nu_\alpha(A) = \int_{\R^n\smallsetminus\set0} \|\bm x\|_2^\alpha \cdot \bm1_A(\bm x/\|\bm x\|_2) \d\mu(\bm x).\end{equation}
$\nu_\alpha(A)$ can be viewed as a type of summary statistic that captures the probability $\mu$ is drawn from a cone specified by a patch of directions $A$, weighted by the $\alpha^{th}$ power of the distance from the origin. Mathematically, it is simply the pushforward of $\|\bm x\|_2^\alpha\d\mu(\bm x)$ on $\R^n\smallsetminus\set 0$ onto $\Sn$ under the map $\bm x\mapsto \bm x/\|\bm x\|_2$. For instance, if $\mu$ has a Lebesgue density $f$ on $\R^n$, then $\nu_\alpha$ has the density 
\begin{equation*}
    w_\alpha(\bm u) = \int_0^\infty r^{n+\alpha-1}f(r \bm u)\d r.
\end{equation*}
Writing $Q=K^\circ \in \mathscr K_\circ^n$ and using the fact that $h_{Q} = \norm{\cdot}_K$, change of measures gives
\begin{equation*}
    \stackrel[X\sim \mu]{}{\E} h_Q(X)^\alpha  = \int_{\R^n}h_Q(\bm x)^\alpha \d\mu(\bm x) = \int_{\Sn} h_Q(\bm u)^\alpha\d\nu_\alpha(\bm u).
\end{equation*}

Assumption \ref{ass:finite-energy}$(i)$ now is equivalent to $\nu_\alpha\left(H_{\bm u}^+\cap \Sn\right) >0 \fa \bm u\in \Sn$. Indeed note that since $\|\bm x\|_2^\alpha >0$ on $H_{\bm u}^+$ we have \[\nu_\alpha(H_{\bm u}^+\cap \Sn)  = 0\iff \int_{H_{\bm u}^+} \|\bm x\|_2^\alpha \d\mu(\bm x) = 0 \iff \mu(H_{\bm u}^+) = 0.\] 
Observing from \eqref{eq:def-nu} that $\nu_\alpha(\Sn) = \E_{X\sim \mu}\|X\|_2^\alpha$ directly shows that Assumption \ref{ass:finite-energy}$(ii)$ is equivalent to $\nu_\alpha(\Sn) < \infty$.

Let us now summarize this discussion. The objective and volume, expressed via the polar, are 
\begin{equation} \label{eq:J-and-Phi}
\mathcal J_\alpha(Q) \sett \int_{\Sn} h_Q^\alpha \d\nu_\alpha, \qquad \qquad \Phi(Q) \sett \vol(Q^\circ) = \frac1n\int_{\Sn}h_Q^{-n}\d\sigma,
\end{equation}
where $\alpha \ge 1$. These are defined over $\mathscr K_\circ^n$.  With an abuse of notation, we shall identify the convex bodies $Q$ via its support function $h=h_Q$ and we shall write $h\in \mathscr K_\circ^n, \cJ_\alpha(h), \Phi(h)$ respectively for $Q\in \mathscr K_\circ^n, \cJ_\alpha(Q), \Phi(Q)$.  Subsequently, \eqref{eq:second-reform} may be formulated as the following
\begin{equation}\label{eq:final_optimization_problem}
\min\ \cJ_\alpha(h) \qquad \text{ s.t. } \qquad h\in \mathscr K_\circ^n, \Phi(h) \le 1.
\end{equation}


Assumption \ref{ass:finite-energy} can be reformulated into the following equivalent assumption on the radial measure $\nu_\alpha$.
\begin{ass}\label{ass:assumption-nu}$\nu_\alpha$ is a positive Borel measure on $\Sn$ satisfying $(i)\ \nu_\alpha\left(H_{\bm u}^+\cap \Sn\right) >0 \fa \bm u\in \Sn$, and $(ii)\ \nu_\alpha(\Sn) < \infty$.
\end{ass}

Impose Assumption \ref{ass:assumption-nu} throughout. The arguments below use only these two properties of $\nu_\alpha$. In particular, $\nu_\alpha$ is allowed to be an atomic measure as long as Assumption \ref{ass:assumption-nu} is satisfied.

\subsubsection{Some preparatory results}  Our next step is to establish a sequence of preparatory results that are useful for showing the existence of flows.  We identify the set of convex bodies as $\mathscr K_\circ^n \subseteq C(\Sn) \subseteq L^2(\Sn,\sigma)$.  Subsequently, the set of feasible solutions may be written as
\begin{equation*}
    \mathscr{A} \sett \set{h\in \body n : \Phi(h)\le 1}.    
\end{equation*}
The optimization instance \eqref{eq:final_optimization_problem} can be expressed as an unconstrained minimization over the energy functional $\cE_\alpha = \cJ_\alpha + \iota_{\mathscr{A}}$ whereby we absorb the constraint as a $0-\infty$ indicator of $\mathscr{A}$; that is, $\cE_\alpha(h)=\cJ_\alpha(h)$ if $h\in \mathscr{A}$, and $\cE_\alpha(h)=\infty$ otherwise. Our next two steps are to verify Assumption \ref{ass:needed-assumptions} for the choices $\cE=\cE_\alpha$ and $S=\mathscr{A}$.  It is clear that $\cJ_\alpha$ has finite values on $\mathscr{A}\ne \varnothing$ from which it follows that $\cE_\alpha$ is proper.  Next, we show that $\cE_\alpha$ is also lower semicontinuous. 
\begin{lem}\label{lem:lsc}
$\cE_\alpha$ is lower semicontinuous in $L^2(\Sn,\sigma)$.
\end{lem}
\begin{proof}
Suppose $h_j\in \mathscr{A}$ and $h_j\to h$ in $L^2$. Proposition \ref{prop:equivalent-metrics} gives a subsequence that converges uniformly to a nonnegative support function representing $h$ in $L^2(\Sn,\sigma)$. The circumradii of the elements of this subsequence are bounded since $\|h_Q\|_\infty = \max_{\bm x\in Q}\|\bm x\|_2$. Say all these circumradii are bounded above by $R>0$. This gives that all these support functions are $R-$Lipschitz. Indeed say $\|h_Q\|_\infty \le R$ and pick $\bm u,\bm v\in \Sn$ with $h_Q(\bm u) \ge h_Q(\bm v)$ and pick $\bm x_{\bm u}\in Q$ with $h_Q(\bm u) = \inner{\bm x_{\bm u}}{\bm u}$. So $h_Q(\bm v) \ge \inner{\bm x_{\bm u}}{\bm v}$. Then $h_Q(\bm u) - h_Q(\bm v) \le \inner{\bm x_{\bm u}}{\bm u-\bm v} \le \|\bm x_{\bm u}\| \|\bm u-\bm v\|_2 \le R \|\bm u-\bm v\|_2$. Therefore a common circumradius gives a common Lipschitz constant.

Fatou gives
\begin{equation}\label{eq:fatou-volume}
    \frac1n\int h^{-n}\d\sigma\le1.    
\end{equation}
If $h(\bm u_0)=0$, then $h(\bm u)\le R\|\bm u-\bm u_0\|$. A cap of radius $\epsilon$ around $\bm u_0$ gives
\begin{equation*}
    \int h^{-n}\d\sigma \ge c_n\epsilon^{n-1}(R\epsilon)^{-n} = \frac{c_n}{R^n\epsilon}.
\end{equation*} 
Taking $\epsilon\downarrow0$ contradicts \eqref{eq:fatou-volume}. Thus $h>0$ and $h\in \mathscr{A}$. By Proposition \ref{prop:equivalent-metrics}, the original sequence $(h_j)$ converges uniformly to $h$, since it converges in $L^2$ and all these functions are support functions which are unique in $L^2(\Sn)$ due to $\sigma$ having full support. Hence $\mathscr{A}$ is closed in $H$. Uniform convergence and finiteness of $\nu_\alpha$ give $\abs{\cJ_\alpha(h_j) - \cJ_\alpha(h)} \le \nu_\alpha(\Sn) \|h_j^\alpha-h^\alpha\|_\infty \to 0$. Since $\mathscr{A}$ is closed in $H$, we conclude the lower semicontinuity of $\cE_\alpha$.
\end{proof}

Next we know that $\mathscr{A}$ is convex by Proposition \ref{prop:constr-convex} and verify that $\cJ_\alpha$ is convex on $\mathscr{A}$. This result automatically shows $\cE_\alpha$ is convex on $H$. Additionally, strict convexity and strong convexity of $\cJ_\alpha$ hold under the additional assumptions stated in the lemma below.

\begin{lem}
\label{lem:polar-convexity}
\begin{enumerate}[label=$(\roman*)$]
\item $\cJ_\alpha$ is convex on $\mathscr{A}$ and strictly convex if $\alpha>1$ and $\operatorname{supp}\nu_\alpha= \Sn$. 
\item Suppose $\alpha=2$ and $\nu_2 \ge w_-\sigma$ for some $w_->0$. Let $h_0,h_1 \in \mathscr{A}$, and $h_t\sett th_0+(1-t)h_1\fa t\in[0,1]$. Then, writing $\cJ\equiv \cJ_2$, we have \begin{equation}
\cJ(h_t) \le t\cJ(h_0) + (1-t) \cJ(h_1) - w_-t(1-t)\norm{h_0-h_1}_2^2.
\label{eq:J2-strong-convexity}
\end{equation}
\end{enumerate}
\end{lem}

\begin{proof}
Fix any $h_0,h_1 \in \mathscr{A}$, and let $h_t\sett th_0+(1-t)h_1\fa t\in[0,1]$. 
\begin{enumerate}[label=$(\roman*)$]
\item \begin{align}
\cJ_\alpha(h_t) &= \int_{\Sn} h_t^\alpha \d\nu_\alpha \notag\\
&\le t\int_{\Sn} h_0^\alpha \d\nu_\alpha + (1-t)\int_{\Sn} h_1^\alpha \d\nu_\alpha\qquad [\because (0,\infty)\ni z\mapsto z^\alpha \text{ is convex for }\alpha\ge 1]\label{eq:convex-s^alpha}\\
&= t\cJ_\alpha(h_0) + (1-t) \cJ_\alpha(h_1).\notag
\end{align}
If $\alpha>1, t\in(0,1)$ and $h_0\ne h_1$, continuity implies that $h_0$ and $h_1$ differ on a nonempty open subset $W$ of $\Sn$. With $\operatorname{supp}\nu_\alpha=\Sn$, $\nu_\alpha(W) > 0$. Strict convexity of $z\mapsto z^\alpha$ makes \eqref{eq:convex-s^alpha} a strict inequality. 
\item Now we take $\alpha=2$ and write $\cJ\equiv \cJ_2$. Write $s=1-t$. From the identity \[(ta+sb)^2 = t^2a^2+s^2b^2 + 2stab = (1-s)ta^2 + (1-t)sb^2 + 2stab = ta^2+sb^2-st(a-b)^2\] we have
\begin{align*}
\cJ(h_t) &= \int_{\Sn} (th_0^2+sh_1^2 - st(h_0-h_1)^2) \d\nu_2 \\
&= t\cJ(h_0) + (1-t)\cJ(h_1) - (1-t)t\int_{\Sn}(h_0-h_1)^2\d\nu_2 \\
&\le t\cJ(h_0) + (1-t)\cJ(h_1) - w_-t(1-t)\norm{h_0-h_1}_2^2.
\end{align*}
\end{enumerate}
\end{proof}

\begin{remark}
The hypothesis in the strict convexity result in Lemma \ref{lem:polar-convexity}$(i)$ above does not hold for atomic measures on finitely many points since it requires full support for measure being evaluated against. They can still satisfy Assumption \ref{ass:assumption-nu} as Section \ref{sec:atomic-exs} shows.
\end{remark}

\begin{lem}\label{lem:power-n-convex}
Let $\lambda_R\sett n(n+1)R^{-n-2}$ and $f(x) = x^{-n}$. For any $x,y\in (0,R],t\in [0,1]$ we have \[f((1-t)x+ty) \le (1-t)f(x) + tf(y) - \frac{\lambda_R}2 t(1-t)(x-y)^2.\]
\end{lem}
\begin{proof}
This immediately follows from $f''(x) = n(n+1)x^{-n-2} \ge \lambda_R > 0$ for $x\in(0,R]$. 
\end{proof}

Define
\[\zeta_{\nu_\alpha} \sett  \min_{\bm x\in S^{n-1}} \int_{\Sn}\inner{\bm u}{\bm x}_+^\alpha \d\nu_\alpha(\bm u), \qquad a_+\sett \max\set{a,0}.\]
The integral is continuous in $\bm x$ by dominated convergence and strictly positive for every $\bm x$ by Assumption \ref{ass:assumption-nu}. Since $\Sn$ is compact, we conclude $\zeta_{\nu_\alpha}>0$.

For a feasible support function $h_Q$ ($Q\in \body n$), choose $R\bm x\in Q$ with $\|\bm x\|_2=1$ where $R=\max_{\bm y\in Q} \|\bm y\| = \|h_Q\|_\infty$. Since $Q\in \body n$, $h_Q(\bm u) \ge \inner{R\bm x}{\bm u}_+ = R \inner{\bm x}{\bm u}_+\fa \bm u\in \Sn$. Hence 
\begin{equation}
\cJ_\alpha(h_Q) \ge \zeta_{\nu_\alpha} \|h_Q\|_\infty^\alpha.\label{eq:energy-lower-bound-support}
\end{equation}
 
\begin{lem}
$\cE = \cE_\alpha$ satisfies Assumption \ref{ass:needed-assumptions} with $S=\mathscr{A}, \Psi=\Phi$.
\end{lem}
\begin{proof}
We verify the four assumptions in Assumption \ref{ass:needed-assumptions} one by one as follows.
\begin{enumerate}[label=$(\roman*)$]
\item Lower semicontinuity is given by Lemma \ref{lem:lsc} and the convexity is supplied by Lemma \ref{lem:polar-convexity}$(i)$. 
\item Note that $h>0$ on $\Sn$ if $h\in \mathscr{A}$. $\nu_\alpha$ is nonzero. Hence 
\begin{equation*}
    \cJ_\alpha(th) = \int_{\Sn} t^\alpha h^\alpha \d\nu_\alpha < \int_{\Sn} h^\alpha \d\nu_\alpha = \cJ_\alpha(h) \qquad \fa t\in(0,1), th\in \mathscr{A}.
\end{equation*}
\item $(0,\infty)\ni z\mapsto z^{-n}$ is strictly convex by Lemma \ref{lem:power-n-convex} whence $\Phi$ is strictly convex by the proof of Lemma \ref{lem:polar-convexity}$(i)$. We further have 
\begin{equation*}
    \Phi(th) = \frac1n\int_{\mathbb S^{n-1}}(th)^{-n}\d\sigma = t^{-n}\cdot \frac1n\int_{\mathbb S^{n-1}}h^{-n}\d\sigma = t^{-n}\Phi(h).
\end{equation*}
Therefore $\mathscr{A} = S = \operatorname{dom}\cE_\alpha = \set{h\in \mathscr K_\circ^n : \Phi(h)\le 1}$. 
\item $\mathscr{A}_c = \set{h\in \mathscr K_\circ^n : \Phi(h)\le 1, \cE_\alpha(h) \le c} = \set{h\in \mathscr K_\circ^n : \Phi(h)\le 1, \cJ_\alpha(h) \le c}$ is $L^p$ bounded for $p=\infty$. Indeed, for every $h\in \mathscr{A}_c$ we have $\|h\|_\infty \le (c/\zeta_{\nu_\alpha})^{1/\alpha}$ by \eqref{eq:energy-lower-bound-support}.
\end{enumerate}
\end{proof}

Therefore the problem \eqref{eq:final_optimization_problem} has a unique minimizer $h_\alpha^* \in \mathscr{A}$ by Lemma \ref{lem:unique-min}. For this particular energy, existence and uniqueness also follow from \cite[Theorem~3.1]{luo2020polar}, applied with $\varphi(t)=t^\alpha$ and measure $\nu_\alpha$.

\subsection{Gradient flow theory for optimal regularization}

\subsubsection{$L^2$-gradient flows for optimal regularization}  We proceed to describe the specialization of $L^2$-gradient flows developed in Section \ref{sec:L2-grad} within the context of optimal regularization.  To do so, we view the set of feasible solutions as $\mathscr{A}\subseteq H \sett (L^2(\Sn), \norm\cdot_2)$ -- here, the crucial detail is that $H$ is a Hilbert space.  By Proposition \ref{prop:iteration-defined-euler}, the discrete Euler scheme
\begin{equation*}
    h_{k+1} = \arg\min_{h\in \mathscr{A}} \set{\cJ_\alpha(h) + \frac{\|h_k-h\|_2^2}{2\tau}}.
\end{equation*}
starting at $h_0\in \mathscr{A}$ is well defined (in the sense that the update exists uniquely), satisfies \eqref{eq:l2-euler-differential-gen}--\eqref{eq:l2-euler-convergence-gen} and exhibits uniform convergence $h_k\to h_\alpha^*$. By Proposition \ref{prop:L2-flow-general}, there is a continuous flow $g:[0,\infty)\to H$ starting at $h_0\in \mathscr{A}$ satisfying the differential inclusion $-\dot g(t) \in \partial \cE_\alpha(g(t))$ for almost every $t>0$ with all $g(t)\in \mathscr{A}$, having regularity \eqref{eq:l2-regularity}, and satisfies the variational inequality \eqref{eq:EVI-convex} and energy dissipation identity \eqref{eq:energy-dissipation}. Further this flow has uniform convergence to the optimal point, namely $\norm{g(t)-h_\alpha^*}_\infty \to 0$ as $t\to\infty$ with energy gap \eqref{eq:flow-sublinear-rate}.  

Consider the special case where $\alpha=2$ and suppose further $\nu_2\ge w_-\sigma$ for some $w_->0$.  Then Lemma \ref{lem:polar-convexity} $(ii)$ gives $\cE_2$ is $2w_--$strongly convex, which then gives a rate of $\norm{g(t) - h_2^*}_2 \le e^{-2w_-t}\norm{h_0-h_2^*}_2$ from \eqref{eq:flow-distance-rate} of Lemma \ref{lem:quadratic-flow-rate}. In order to get rates for general $\alpha$, we need to verify the quadratic growth assumption \eqref{eq:qg}. 

Focus on the sub-level set $\mathscr{A}_{\cJ_\alpha(h_0)}$ where the dynamics starts at $h_0\in \mathscr{A}$. By Corollary \ref{cor:unif-ball}, $\exists 0<r_0\le R_0$ such that $r_0\le h \le R_0\fa h\in \mathscr{A}_{\cJ_\alpha(h_0)}$. We show that the \eqref{eq:qg} is satisfied with 
\begin{equation*}
    \gamma_\alpha \sett \frac{\alpha\cJ_\alpha(h_\alpha^*)(n+1)}{2nR_0^{n+2}}.    
\end{equation*}

\begin{lem}\label{lem:polar-quadratic-growth}
$\cJ_\alpha(h)-\cJ_\alpha(h_\alpha^*) \ge \gamma_\alpha \norm{h-h_\alpha^*}_2^2$ for all $h\in \mathscr{A}_{\cJ_\alpha(h_0)}$.
\end{lem}
\begin{proof}
Throughout this proof, all integrals are over $\Sn$ with respect to the unnormalized surface measure $\sigma$ if the measure is unspecified. Fix $h\in  \mathscr{A}_{\cJ_\alpha(h_0)}$. Since $\cJ_\alpha(h_\alpha^*) \le \cJ_\alpha(h_0)$, both $h$ and $h_\alpha^*$ belong to $\mathscr{A}_{\cJ_\alpha(h_0)}$ and satisfy $r_0\le h,h_\alpha^*\le R_0$. Lemma \ref{lem:unique-min} guarantees that $\Phi(h_\alpha^*)=1$.  Let $h_t\sett (1-t)h_\alpha^* + th$ and $b_t\sett \Phi(h_t)$ for $t\in[0,1]$. Then
\begin{align}
b_t = \frac1n\int h_t^{-n} &\stackrel{\text{Lemma }\ref{lem:power-n-convex}}{\le} \frac{1-t}n\int (h_\alpha^*)^{-n} + \frac tn \int h^{-n} - \frac{\lambda_{R_0}}{2n}t(1-t)\int (h-h_\alpha^*)^2 \notag\\
&\le 1 - \frac{n+1}{2R_0^{n+2}}t(1-t)\int (h-h_\alpha^*)^2\notag\\
&= 1- \frac{n+1}{2R_0^{n+2}}t(1-t)\norm{h-h_\alpha^*}_2^2.\label{eq:bt-bound}
\end{align}
Let $\tilde h_t\sett b_t^{1/n}h_t$. Then $\Phi(\tilde h_t)=1$ whence $\tilde h_t\in \mathscr{A}$. This gives \begin{align*}
\cJ_\alpha(h_\alpha^*) &\le \cJ_\alpha(\tilde h_t) = b_t^{\alpha/n} \cJ_\alpha(h_t)\\
\implies b_t^{-\alpha/n} \cJ_\alpha(h_\alpha^*) &\le \cJ_\alpha(h_t) \stackrel{\text{convexity}}{\le} (1-t)\cJ_\alpha(h_\alpha^*) + t\cJ_\alpha(h) = \cJ_\alpha(h_\alpha^*) + t(\cJ_\alpha(h)-\cJ_\alpha(h_\alpha^*))\\
\implies \cJ_\alpha(h)-\cJ_\alpha(h_\alpha^*) &\ge \frac{b_t^{-\alpha/n}-1}{t} \cJ_\alpha(h_\alpha^*) \\
&\ge \frac{-\alpha(b_t-1)}{nt}\cJ_\alpha(h_\alpha^*) \qquad\qquad [\because x^{-p}-1\ge -p(x-1) \fa p>0,x>0]\\
&\stackrel{\eqref{eq:bt-bound}}{\ge} \frac{\alpha(n+1)\cJ_\alpha(h_\alpha^*)}{2nR_0^{n+2}} (1-t)\norm{h-h_\alpha^*}_2^2 \\
&= (1-t)\gamma_\alpha\norm{h-h_\alpha^*}_2^2
\end{align*}
Since there were some steps with division by $t$, the above steps hold for every $t\in (0,1]$. Taking limits $t\to 0$ proves the lemma.
\end{proof}

Therefore, the exponential convergence rates \eqref{eq:rate-euler-gen}, \eqref{eq:rate-flow-gen} hold with $\gamma=\gamma_\alpha$ respectively for the discrete scheme and the flow, for any $\alpha\ge 1$. 


\subsubsection{Hausdorff minimizing movement scheme for optimal regularization}  Second, we describe the minimizing movement scheme developed in Section \ref{sec:hausdorff-general} within the context of optimal regularization. In this instance, we view the set of feasible solutions as $\mathscr{A}\subseteq B=(C(\Sn),\norm{\cdot}_\infty)$ to obtain the following Hausdorff minimizing movement scheme.  
\begin{equation}\label{eq:polar-hausdorff-update}
h_{k+1}\in \arg\min_{h\in\mathscr{A}} \set{\cJ_\alpha(h)+\frac{\norm{h_k-h}_\infty^2}{2\tau}}, \qquad h_0\in\mathscr{A}.
\end{equation}
We note that while $\delta_2$ and $\delta_{\infty}$ induce the same topology on the space of support functions, their proximal penalties measure deviations differently and can produce different updates. By Proposition~\ref{prop:euler-scheme-converge}, every update \eqref{eq:polar-hausdorff-update} exists and every choice of an exact update converges uniformly to $h_\alpha^*$. On the sub-level set $\mathscr{A}_{\cJ_\alpha(h_0)}$, Lemma~\ref{lem:polar-quadratic-growth} gives
\[\cJ_\alpha(h)-\cJ_\alpha(h_\alpha^*) \ge\gamma_\alpha\norm{h-h_\alpha^*}_2^2.\]
Let $R_0$ be the common circumradius bound for this sub-level set from Corollary \ref{cor:unif-ball}, and write $\kappa_{n-1}=\vol_{n-1}(B_2^{n-1})$. \cite[Theorem 6]{stability-convex-body}  gives
\[\norm{h-h_\alpha^*}_2^2 \ge\frac{2\kappa_{n-1}(2R_0)^{1-n}}{n(n+1)}\norm{h-h_\alpha^*}_\infty^{n+1}.\]
Thus \eqref{eq:hausdorff-power-growth} holds with $p=n+1$ and 
\begin{equation}\label{eq:polar-hausdorff-growth-constant}
c_\alpha=\gamma_\alpha\frac{2\kappa_{n-1}(2R_0)^{1-n}}{n(n+1)}.
\end{equation}
By taking the limit $\tau\downarrow 0$, one can also find a locally absolutely continuous curve $g:[0,\infty)\to B$ with $g(0) = h_0$ by Proposition \ref{prop:hausdorff-scheme-limit}, which can be approximated via the interpolations \eqref{eq:hausdorff-interpolate}.

\begin{cor}\label{cor:polar-hausdorff-rates}
Let $n\ge 2$. For every fixed $\tau>0$, the Hausdorff iterates for $\cJ_\alpha$ satisfy
\begin{equation}\label{eq:polar-hausdorff-rates}
\begin{aligned}
\cJ_\alpha(h_k)-\cJ_\alpha(h_\alpha^*) &=\mathcal O\left((1+k\tau)^{-(n+1)/(n-1)}\right),\\
\norm{h_k-h_\alpha^*}_\infty &=\mathcal O\left((1+k\tau)^{-1/(n-1)}\right).
\end{aligned}
\end{equation}
Every limit curve satisfies the corresponding bounds with $h_k$ replaced by $g(t)$ and $k\tau$ by $t$.
\end{cor}

\begin{proof}
We apply Proposition~\ref{prop:hausdorff-growth-rates} with $p=n+1,c=c_\alpha$. Then $1/\rho=(n+1)/(n-1)$ and $1/(p\rho)=1/(n-1)$.
\end{proof}

The finite-slope bound also  gives that the resulting curve $g$ is $L_0-$Lipschitz by \eqref{eq:hausdorff-finite-slope} with the bound 
\begin{equation}\label{eq:polar-hausdorff-initial-slope}
L_0 = s(h_0)\le \alpha\norm{h_0}_\infty^{\alpha-1}\int_{\Sn} \d\nu_\alpha(u) = \alpha\|h_0\|_\infty^{\alpha-1}\nu_\alpha(\Sn).
\end{equation}
obtained from \eqref{eq:slope-def} and the gradient inequality $b^\alpha-a^\alpha\le \alpha b^{\alpha-1}(b-a)\fa a,b\ge 0$ for $x\mapsto x^\alpha, x\in \R_{\ge 0}$.

Corollary~\ref{cor:polar-hausdorff-rates} gives the refined energy-gap and Hausdorff-distance rates \eqref{eq:polar-hausdorff-rates}. For the continuous-time limit, use \eqref{eq:hausdorff-interpolate}. Proposition~\ref{prop:hausdorff-scheme-limit} gives a subsequence of every vanishing sequence of time steps for which \eqref{eq:hausdorff-subsequence-limit} holds. Every such limit $g$ is a curve of maximal slope, satisfies \eqref{eq:hausdoff-energy} and \eqref{eq:hausdorff-differential}, and obeys the continuous bounds in Corollary~\ref{cor:polar-hausdorff-rates}.

%% file: sec_Examples.tex
\section{Examples of Optimal Convex Regularizers} \label{sec:optimal-convex-reg-examples}
This section aims to derive explicit, closed-form expressions for the optimal convex regularizer using the KKT theory developed in Section \ref{sec:optimality}. For $\alpha \geq 1$, we consider characterizing the solution to \begin{align}\label{eq:opt-problem}
\min_{K \in \body n} \stackrel[X\sim \mu]{}{\E}\|X\|_K^\alpha \qquad \text{ s.t. } \qquad \vol(K) \leq 1.
\end{align}
We consider distributions $\mu$ that satisfy Assumption \ref{ass:finite-energy}, which allows us to reformulate the problem equivalently as \begin{align}\label{eq:opt-problem-polar}
\min_{Q \in \body n} \stackrel[X\sim \mu]{}{\E}h_Q(X)^\alpha \qquad \text{ s.t. } \qquad \Phi(Q) \leq 1.
\end{align}

Our first example addresses an open problem raised in \cite[\S 8]{leong2025optimal}. When the optimal regularizer for a data distribution has a nonconvex unit ball, which convex regularizer is the best surrogate for the same data? For a density proportional to $\exp(-\norm{\bm x}_L)$, \cite[\S 3.5]{leong2025optimal} shows that the optimal regularizer among gauges of star bodies is a rescaling of $\|\cdot\|_L$. When $L$ is nonconvex, its convex hull suggests a candidate for the optimal convex unit ball. In which cases this candidate indeed minimizes the objective over all convex bodies is unknown, and the convex-constrained problem has only been studied computationally in two-dimensional examples in \cite[\S 4]{leong2025optimal-uncertainty}. 

Section \ref{sec:full-supp-ex} resolves this question for the family of densities proportional to $\exp(-\norm{\bm x}_p)$ for $p\in (0,1)$. Here the optimal star body gauge is proportional to the nonconvex $\ell_p$-quasinorm, and its unit ball has the $\ell_1$-ball as its convex hull, up to scale. Theorem \ref{thm:ell1-ball-opt} shows that for every $\alpha \geq 1$, the volume-normalized $\ell_1$-ball is the unique minimizer of \eqref{eq:opt-problem}. 

Section \ref{sec:atomic-exs} considers data supported on a finite number of atoms given by signed basis elements. In this case, we solve the optimality conditions explicitly and show that the optimal polar body is a shifted zonotope whose edge lengths and translations are determined by the atom weights. When the atoms are the signed coordinate vectors and all have equal weight, the corresponding primal unit ball is again the volume-normalized $\ell_1$-ball.

\subsection{From the $\ell_p$ quasinorm to optimal $\ell_1$ regularization}\label{sec:full-supp-ex}
Let $p\in(0,1)$ and $\mu$ be the measure on $\R^n$ with density $f(\bm x) = Z^{-1} \exp(-\|\bm x\|_p)$ and normalizing constant \[Z \sett \int_{\R^n} \exp(-\|\bm x\|_p)\d\bm x.\] Note that Assumption \ref{ass:assumption-nu} is satisfied because $\operatorname{supp}\mu = \R^n$ and the $L^2$ norm has finite $\alpha-$moment. Hence the reduction from $\mathscr K^n$ to $\body n$ and then converting to polar variables, as in \eqref{eq:opt-problem-polar}, is valid.

\paragraph{\bf Proof roadmap} We prove Theorem \ref{thm:ell1-ball-opt} by verifying the KKT conditions for the volume-normalized $\ell_1$-ball. Note that each facet of the $\ell_1$-ball is a signed copy of the probability simplex, so we can express the first variations of the objective and the volume constraint as averages on the simplex. The data distribution gives a weighted average, while the volume constraint gives the uniform average. The main step is to verify dual feasibility by showing that the weighted average is at least as large as the uniform average for every support function. Lemmas \ref{lem:positive-cov} and \ref{lem:uniform-stable} establish this comparison using a covariance inequality and the permutation symmetry of the weight. Since the $\ell_1$-norm is constant on each facet, the two averages agree when evaluated there, completing the optimality certificate.  

Let $\Delta_{n-1}$ denote the probability simplex \[\Delta_{n-1} \sett \set{\bm x \in \R^n \st \bm x\ge 0, \bm 1^\top \bm x=1}.\]
Denote by $\eta_\Delta$ the uniform probability measure on $\Delta_{n-1}$, $\bm \epsilon \in \set{\pm 1}^n$ the uniform Rademacher, and $\eta_w$ the probability measure on $\Delta_{n-1}$ with $\d\eta_w(\bm s) = D^{-1} \|\bm s\|_p^{-(n+\alpha)}\d\eta_\Delta(\bm s)$ with normalizing constant \[D \sett \int_{\Delta_{n-1}} \|\bm s\|_p^{-(n+\alpha)}\d \eta_\Delta(\bm s).\]

\begin{lem}\label{le:reduction-p-norm}
With $n,p,\alpha,\mu$ as above and $Q\in \body n$, we have \[\stackrel[X\sim \mu]{}{\E}h_{Q}(X)^\alpha = \stackrel[X\sim\mu]{}{\E}\|X\|_1^\alpha \cdot \stackrel[\bm \epsilon]{}{\E}\stackrel[\bm s\sim \eta_w]{}{\E} h_{Q}(\bm \epsilon\odot\bm s)^\alpha.\]
\end{lem}
\begin{proof}
Write $K=Q^\circ, h=h_{Q}, \bm x=r(\bm\epsilon\odot\bm s)$ with  where $r=\|\bm x\|_1$, $\bm \epsilon \in \set{\pm1}^n$ is the sign pattern of $\bm x$, $\bm s=\bm \epsilon\odot \bm x/r \in \Delta_{n-1}$ and $\odot$ represents coordinate-wise multiplication. Under $\mu$, the sign vector is uniform and independent of $(r,\bm s)$. The simplex variable $s$ has density proportional to $\|\bm s\|_p^{-n}$ with respect to $\eta_\Delta$. Under the law obtained by reweighting $\mu$ by $\|X\|_1^\alpha/\mathbb E_\mu\|X\|_1^\alpha$, the simplex variable has law $\eta_w$. Then we have
\begin{align*}
\stackrel[X\sim\mu]{}{\E}\|X\|_K^\alpha &= \frac1Z \int_{\R^n\smallsetminus\set0} h(\bm x)^\alpha e^{-\|\bm x\|_p} \d\bm x \\
&= \frac1Z\sum_{\bm \epsilon\in\set{\pm1}^n}\int_{\bm s\in \Delta_{n-1}} \int_0^\infty h(\bm \epsilon\odot \bm s)^\alpha r^{n-1+\alpha} e^{-r\|\bm s\|_p}\d r \d s_1\cdots \d s_{n-1}\\
&= \frac1Z\sum_{\bm \epsilon\in\set{\pm1}^n}\int_{\bm s\in \Delta_{n-1}}h(\bm \epsilon\odot \bm s)^\alpha \left(\int_0^\infty r^{n-1+\alpha} e^{-r\|\bm s\|_p}\d r \right) \d s_1\cdots \d s_{n-1}\\
&= \frac{2^n\Gamma(n+\alpha)}{(n-1)!Z}\int_{\bm s\in \Delta_{n-1}} \|\bm s\|_p^{-(n+\alpha)}\E_{\bm \epsilon}h(\bm \epsilon\odot \bm s)^\alpha  \d\eta_{\Delta}(\bm s)\\
&= \frac{2^n\Gamma(n+\alpha)D}{(n-1)!Z}\E_{\bm s\sim \eta_w}\E_{\bm \epsilon} h(\bm \epsilon\odot\bm s)^\alpha.
\end{align*}
In the above we used that the change of measures from $(\R_{\ge 0})^n\smallsetminus\set0$ to $\R_{>0}\times \Delta_{n-1}$ has the Jacobian determinant $r^{n-1}$ and $\d\eta_\Delta(\bm s) = (n-1)!\d s_1\cdots\d s_{n-1}$. Note that taking $Q = (\mathcal B_1^n)^\circ = \mathcal B_\infty^n$ so that $h_{Q} = \|\cdot\|_1$ in the above gives the preceding constant $\E_{X\sim \mu} \|X\|_1^\alpha = 2^n\Gamma(n+\alpha)D/((n-1)!Z)$.
\end{proof}

For our analysis, we introduce a class of translated and stretched simplices in $\R^n$ parameterized by $a,b\in \R^+$ with $b>na$: 
\[T_n^{a,b}\sett \set{\bm s\in \R^n : s_i\ge a\fa i, \bm 1^\top\bm s=b} = a\bm 1 + (b-an)\Delta_{n-1}.\] 
Taking $a=0,b=1$ makes $T_{n}^{0,1} = \Delta_{n-1}$. There is a natural permutation action $\fS_n\curvearrowright T_n^{a,b}$ of the symmetric group $\fS_n$ on the simplex, namely $\pi \bm s = (s_{\pi^{-1}(1)}, \cdots, s_{\pi^{-1}(n)})$ for $\pi\in \fS_n, \bm s\in T_n^{a,b}$. In this context, we say that a function $f:T_n^{a,b} \to \R$ is symmetric if $f(\pi\bm s) = f(\bm s)$ for all $\pi\in \fS_n, \bm s\in T_n^{a,b}$. We move on to show that the covariance of any two continuous convex symmetric functions on a uniformly chosen point of any simplex has non-negative covariance. Our proof is motivated from the techniques used in \cite[Corollary 2.6]{schur-convex}, \cite[Chapter 3 Proposition C.2]{majorization} and \cite[Definition 1.1, Property (P3), Theorem 2.1]{association}.

\begin{lem}\label{lem:positive-cov}
Let $n\in \N, a,b\in \R$ be such that $b> na$, and consider the simplex $T_n^{a,b}$. Suppose $F,G:T_n^{a,b}\to \R$ are continuous, convex, and symmetric. Then $\E [F(X)G(X)] \ge \E F(X) \E G(X)$ where $X\sim \operatorname{Unif}(T_n^{a,b})$.
\end{lem}
\begin{proof}
We induct on $n\ge 1$. For $n=1,T_n^{a,b}=\set{b}$. Then $\E [F(X)G(X)] - \E F(X) \E G(X) = 0$.

Now suppose the statement is true for $n=k-1$ and we prove the statement for $n=k$. Let $X\sim \operatorname{Unif}(T_n^{a,b})$. Consider the minimum entry along with its index namely $M=\min X_i, I=\arg\min X_i$. $I$ is unique almost surely. By the law of total expectation, $\E[H(X)] = \E[\E[H(X) \mid M, I]]$ where $H\sett FG$. The inner conditional expectation is taken on the slice inside $T_n^{a,b}$ where $X_i=m$ and $X_j\ge m\fa j\ne i$. More precisely, this is the slice \[T_n^{a,b}(i,m) \sett \set{\bm s\in \R^n : s_i = m, s_j\ge m\fa j\ne i, \sum_j s_j  = b} \cong \set{m} \times T_{n-1}^{m,b-m}\] where the latter is understood as those vectors with $m$ in the $i^{th}$ coordinate. Since $X_i$ is fixed on this slice, we only integrate over $T_{n-1}^{m,b-m}$. Moreover, $\E[\E[H(X) \mid M, I=i]] = \E[\E[H(X) \mid M, I=1]]\fa i$ due to the permutation invariance of $H$, namely $\mathbb P[I=i\mid M] = n^{-1}\fa i\in[n]$. By induction hypothesis, $\E[F(X)G(X) \mid M=m,I=1] \ge \E[F(X) \mid M = m,I=1]\E[G(X)\mid M=m,I=1]$ for almost every $m$, since conditional on $M=m,I=1$, $(X_2,\cdots,X_n)$ is uniform on $T_{n-1}^{m,b-m}$. Let $A(m) = \E[F(X) \mid M=m, I=1], B(m) = \E[G(X) \mid M=m, I=1]$. Therefore,
\begin{equation}
\begin{aligned}
\E[F(X)G(X)] &= \E[\E[F(X)G(X) \st M, I=1]]\\
&\ge \E[\E[F(X) \mid M, I=1]\cdot \E[G(X) \mid M, I=1]]\\
&= \E[A(M) B(M)].
\end{aligned}\label{eq:covariance-solve}
\end{equation}

Now we claim that if $a\le m < m' \le b/n$ then $A(m) \ge A(m'), B(m) \ge B(m')$. We just show for $A$ since $F,G$ have the same properties. First, the minima of $F,G$ occurs at the center $\bm 1b/n$. This is because for any $\bm s\in T$ (same argument for $G$), \[F(\bm 1b/n) = F\left(\frac1{n!}\sum_{\pi\in \fS_n} \pi\bm s\right) \stackrel{\text{convex}}\le \frac1{n!}\sum_{\pi\in \fS_n} F(\pi\bm s) = \frac1{n!}\sum_{\pi\in \fS_n} F(\bm s) \stackrel{\text{symmetry}}= F(\bm s).\] 
\begin{figure}[h]\centering
\includegraphics[trim={0 0 0 1cm},clip,width=9cm]{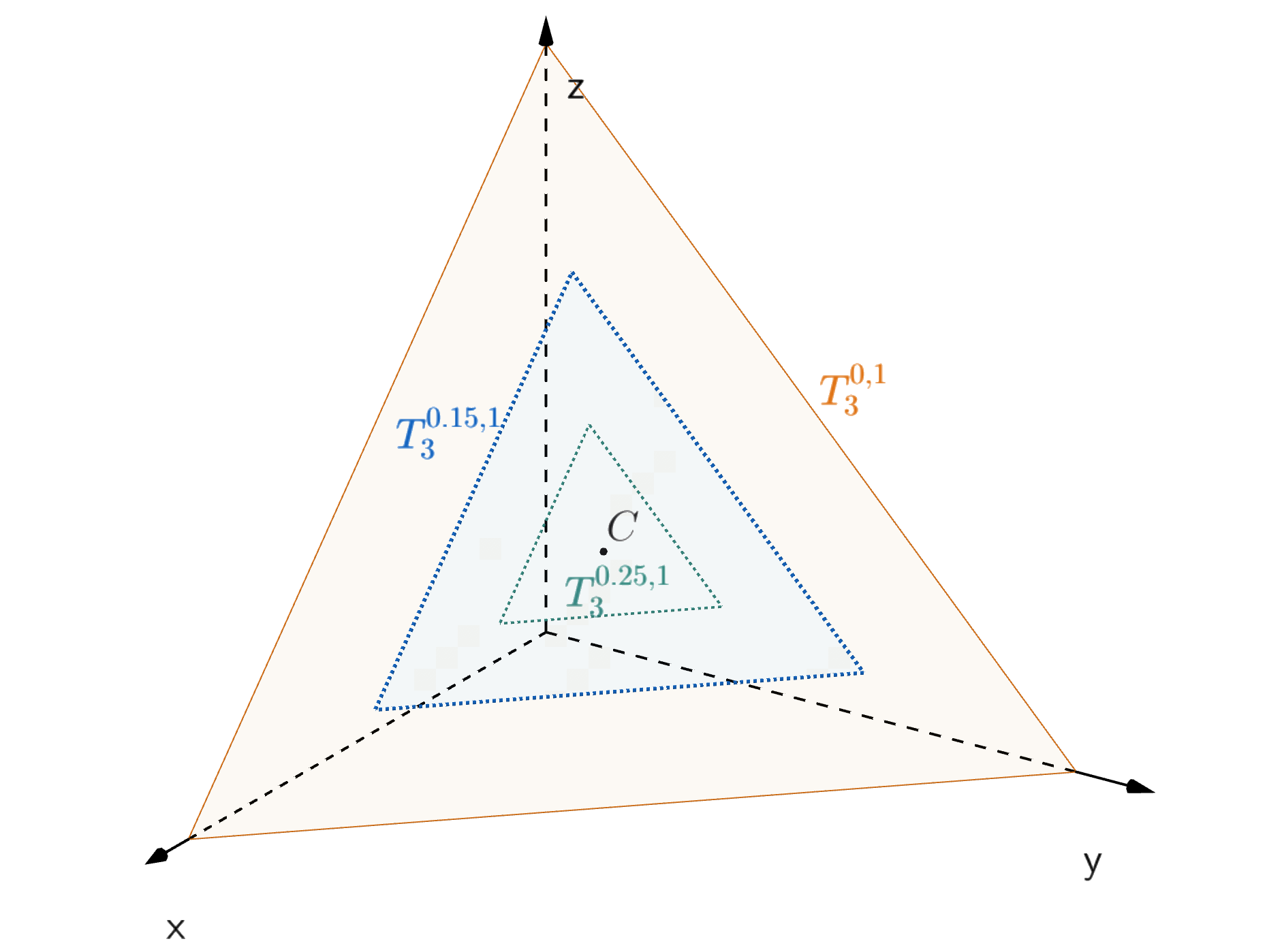}
\caption{Simplices $T_3^{0,1}, T_3^{0.15,1}, T_3^{0.25,1}$. $T_3^{0,1}(1,0.15), T_3^{0,1}(1,0.25)$ are slices of $T_3^{0.15,1}, T_3^{0.25,1}$.}\label{fig:simplices}
\end{figure}

Next refer to Figure \ref{fig:simplices} for better understanding of the following argument. Note that varying $m\in [a,b/n]$, all $T_{n}^{a,b}(1,m)$ are homothetic about the center $\bm 1 b/n$. This is clear from the geometry of $T_n^{a,b}(1,m)$: $T_n^{a,b}(1,m) = m\bm 1 + (b-nm) (\set 0\times \Delta_{n-2})$.
Indeed, the homothety is $H_t(\bm s) = (1-t)\bm 1 b/n + t\bm s$. Fix $a\le m<m'< b/n$. Then for $\theta=(b/n-m')/(b/n-m)\in(0,1)$, the restricted map $H_\theta: T_{n}^{a,b}(1,m) \to T_{n}^{a,b}(1,m')$ is a bijection between the slices. Let's verify that this is well-defined. Any $s\in T_{n}^{a,b}(1,m)$ satisfies $s_j\ge m\fa j>1, \bm 1^\top \bm s = b$ whence $\bm s' = H_\theta(\bm s)$ satisfies that $s'_j = b/n + (s_j-b/n) \theta \ge b/n + (m-b/n)(b/n-m')/(b/n-m) = b/n + m'-b/n = m'$ and $\bm 1^\top \bm s' = b(1-\theta) + \theta b = b$. Bijection follows clearly from the linear interpolation between two distinct points $\bm 1b/n$ and $\bm s$. Therefore, $F(H_\theta(\bm s)) \le (1-\theta)F(\bm 1b/n) + \theta F(\bm s) \le (1-\theta)F(\bm s) + \theta F(\bm s) = F(\bm s)$. Since $H_\theta$ is a linear interpolation, $H_\theta(\bm s)$ is uniform over $T_n^{a,b}(1,m')$ if $\bm s$ is uniform over $T_n^{a,b}(1,m)$. Integrating $F(H_\theta(\bm s)) \le F(\bm s)$ over $T_{n}^{a,b}(1,m)$ gives $A(m') \le A(m)$. Doing the same for $G$ gives $B(m') \le B(m)$. Therefore $A,B$ are non-increasing on $(a,b/n)$. This shows that if $M,\tilde M$ are independent copies of $\min X_i$ where $X\sim \operatorname{Unif}(T)$, then $(A(M)-A(\tilde M))(B(M)-B(\tilde M)) \ge 0$ since the factors are either both non-positive or both non-negative depending on whether $M\ge \tilde M$ or $\tilde M\ge M$, respectively. Therefore $\E_{M,\tilde M}[(A(M)-A(\tilde M))(B(M)-B(\tilde M))] \ge 0$. But this means $\E[A(M)B(M)] \ge \E[A(M)]\E[B(M)]$. From \eqref{eq:covariance-solve}, $\E[F(X)G(X)] \ge \E[A(M)] \E[B(M)] = \E[F(X)]\E[G(X)]$.
\end{proof}

An immediate consequence of the above lemma is that if $F:\Delta_{n-1}\to \R$ is continuous, symmetric and convex then noting that the density $G(\bm s) \sett D^{-1}\|\bm s\|_p^{-(n+\alpha)}$ of $\eta_w$ is also continuous, symmetric and convex on $\Delta_{n-1}$ we get $\E_{X\sim \eta_w} F(X) \ge \E_{\bm s\sim \eta_\Delta} F(X)$.

\begin{lem}\label{lem:uniform-stable}
Let $n,\alpha,p, \eta_\Delta, \eta_w$ be as above and $f:\Delta_{n-1} \to \R$ be continuous and convex. Then \[\E_{\bm s\sim \eta_w} f(\bm s) \ge \E_{\bm s\sim \eta_\Delta} f(\bm  s).\]
\end{lem}
\begin{proof}
We first symmetrize \[F(\bm s) \sett \frac{1}{n!}\sum_{\pi\in \mathfrak S_n} f(\pi \bm s).\] 
$\eta_\Delta,\eta_w$ are $\fS_n-$invariant whence $\E_{\bm s\sim \eta_w} F(\bm s) = \E_{\bm s\sim \eta_w} f(\bm s)$ and $\E_{\bm s\sim \eta_\Delta} F(\bm s) = \E_{\bm s\sim \eta_\Delta} f(\bm s)$. Thus, it suffices to show $\E_{\bm s\sim \eta_w} F(\bm s)\ge \E_{\bm s\sim \eta_\Delta} F(\bm s)$. Note that $F$ is continuous, convex and symmetric. The density $G(\bm s) = D^{-1}\|\bm s\|_p^{-(n+\alpha)}$  of $\eta_w$ is clearly continuous and symmetric. Since $p\in(0,1)$, $\bm s\mapsto \sum_i s_i^p$ is concave and positive, and $r\mapsto r^{-(n+\alpha)/p}$ is convex and decreasing whence their composition $G$ is convex on $\Delta_{n-1}$. By Lemma \ref{lem:positive-cov}, $\E[F(\bm s)G(\bm s)] \ge \E[F(\bm s)] \E[G(\bm s)]$ so $\E_{\bm s\sim \eta_w} F(\bm s)\ge \E_{\bm s\sim \eta_\Delta} F(\bm s)$.
\end{proof}

\begin{thm}\label{thm:ell1-ball-opt}
With $n,p,\alpha,\mu,\eta_w,\eta_\Delta$ as above, the unique minimizer of \eqref{eq:opt-problem} is  $K^* = (n!/2^n)^{1/n} B_{\ell_1}^n$.
\end{thm}
\begin{proof}
Every spherical function is extended by positive homogeneity in this proof. Let $M\sett \E_{X\sim\mu}\|X\|_1^\alpha$. We work with the following polar form \[\min_{Q\in \body n} \stackrel[X\sim \mu]{}{\E}\|X\|_{Q^\circ}^\alpha \qquad \text{s.t.} \qquad \vol(Q^\circ) \le 1.\] By Lemma \ref{le:reduction-p-norm}, the objective $E(h_Q) = E(Q)$ given by \[ E(Q) \sett \stackrel[X\sim \mu]{}{\E}\|X\|_{Q^\circ}^\alpha = M \stackrel[\bm \epsilon]{}{\E}\stackrel[\bm s\sim \eta_w]{}{\E} h_Q(\bm \epsilon\odot\bm s)^\alpha.\]  
By the same argument in Lemma \ref{le:reduction-p-norm}, the constraint function is $\Phi(h) = V_n\cdot \E_{\bm \epsilon}\E_{\bm s\sim\eta_\Delta}h(\bm \epsilon\odot \bm s)^{-n}$ where $V_n = 2^n/n!$. Let $h_*(\bm x) = V_n^{1/n}\|\bm x\|_1$. This is the gauge of the volume normalized $\cB_1^n$.

By Lemma \ref{lem:derivative-of-integral}, their respective derivatives are \[
DE(h)[g] = \alpha M \stackrel[\bm s\sim \eta_w, \bm \epsilon]{}{\E} [h(\bm \epsilon\odot \bm s)^{\alpha-1} g(\bm \epsilon\odot \bm s)], \quad 
D\Phi(h)[g] = -nV_n \stackrel[\bm s\sim \eta_\Delta,\bm \epsilon]{}{\E} [h(\bm \epsilon\odot \bm s)^{-n-1}g(\bm \epsilon\odot \bm s)].\]
Therefore 
\[
DE(h_*)[g] = \alpha MV_n^{(\alpha-1)/n}  \stackrel[\bm s\sim \eta_w,\bm \epsilon]{}{\E} g(\bm \epsilon\odot \bm s), \qquad 
D\Phi(h_*)[g] = -nV_n^{-1/n} \stackrel[\bm s\sim \eta_\Delta,\bm \epsilon]{}{\E} g(\bm \epsilon\odot \bm s).
\]
Let's invoke the KKT conditions from Proposition \ref{prop:kktconditions} applied to the vector space of continuous functions $C(\Sn)$ with the open convex domain $\Omega = \set{f\in C(\Sn) : \min f>0}$, $\cF_0=E+\iota_{\mathscr F^n}, \cF_1=\Phi-1$. First note that Slater's condition holds for the constant support function $h_0\equiv R$ (on $\Sn$) where $R > \vol(\mathcal B_2^n)^{1/n}$ whence $E(h_0) < \infty$ and $\Phi(h_0) = \vol(R^{-1}\cB_2^n) = R^{-n}\vol(\mathcal B_2^n) < 1$. Take $\nu\sett DE(h_*) + \lambda_* D\Phi(h_*)$ with $\lambda_* = \alpha V_n^{\alpha/n}M/n > 0$. Fix any support function $g\in \mathscr F^n$. For each $\epsilon \in \set{\pm 1}^n$, the map $\Delta_{n-1}\to \R$ given by $\bm s\mapsto g(\bm\epsilon\odot \bm s)$ is continuous and convex because $g$ is convex whence Lemma \ref{lem:uniform-stable} gives $\E_{\bm s\sim \eta_w} g(\bm\epsilon \odot \bm s) \ge \E_{\bm s\sim \eta_\Delta} g(\bm\epsilon \odot\bm s)$. Therefore \[\frac{\nu[g]}{\alpha MV_n^{\frac{\alpha-1}{n}}} =  \stackrel[\bm \epsilon]{}{\E}\left(\stackrel[\bm s\sim \eta_w]{}{\E} g(\bm \epsilon\odot\bm s)-\stackrel[\bm s\sim \eta_\Delta]{}{\E} g(\bm \epsilon\odot\bm s)\right) \ge 0.\] 
Hence $\nu$ is dual feasible. Further $\nu[h_*] = 0$ and thus $-\nu\in N_{\mathscr F^n}(h_*)$. Stationarity follows from $DE(h_*)+\lambda_*D\Phi(h_*)=\nu$. Primal feasibility holds because $\Phi(h_*)=1$. Conclude by Proposition \ref{prop:kktconditions}. The strict convexity of $\lambda_*\Phi$ (since $\lambda_*>0$) and convexity $E$ give a unique minimizer.
\end{proof}

\subsection{Shifted zonotopes for atomic distributions} \label{sec:atomic-exs}

We next consider atomic measures. Let $\mu$ be a probability measure on $\R^n$ supported on finitely many nonzero points, say \[\mu = \sum_i a_i \delta_{\bm u_i}\] for some points $P=\set{\bm u_1,\cdots,\bm u_m} \subseteq \R^n\smallsetminus \set0$ where the $a_i$'s are positive and sum to one. Then, \[\stackrel[X\sim \mu]{}{\E}[h_Q(X)^\alpha] = \sum_{i=1}^m a_i h_Q(\bm u_i)^\alpha.\]
The corresponding measure $\nu_\alpha$ from \eqref{eq:def-nu} is \[\nu_{\alpha,P}(A) = \sum_{i=1}^m a_i \|\bm u_i\|_2^\alpha \cdot \bm1_A(\bm u_i/\|\bm u_i\|_2), \qquad A\stackrel[\text{Borel}]{}\subseteq \Sn.\]
This measure satisfies Assumption \ref{ass:assumption-nu} as long as $P$ is not contained in any closed half space through the origin, equivalently $0 \in \operatorname{int}\conv P$. Therefore the gradient flow results of Section \ref{sec:application-to-our-problem} apply to these atomic measures. In particular, owing to the fact that each sublevel set supplies a common circumradius for the support functions it contains, we get a positive strong convexity modulus for $\Phi$ and hence the quadratic growth rates of Lemma \ref{lem:polar-quadratic-growth} apply to the $L^2$ flow yielding exponential convergence rates from Lemma \ref{lem:qg-estimates} for the Euler scheme and the polynomial convergence rates in Corollary \ref{cor:polar-hausdorff-rates} for the Hausdorff scheme.

\begin{remark}
If $\operatorname{span} P \ne \R^n$ then $\exists\ \bm u\in \Sn$ with $\inner{\bm u}{\bm v} = 0\fa \bm v\in P$ so Assumption \ref{ass:assumption-nu}$(i)$ fails. But Assumption \ref{ass:assumption-nu}$(ii)$ is satisfied for any finite $P$ since $\nu_{\alpha,P}(\Sn) = \sum_ia_i\|\bm u_i\|_2^\alpha$.
\end{remark}

For the rest of this section, we take $P=\set{\pm \bm u_1,\cdots, \pm\bm u_n}$ where $\set{\bm u_i}_{i=1}^n\subset \Sn$ is a basis of $\R^n$ and $\mu = \sum_i(a_i\delta_{\bm u_i} + b_i\delta_{-\bm u_i})$ where $a_i,b_i>0$ and sum to one. Clearly $\nu_{\alpha,P}$ satisfies Assumption \ref{ass:assumption-nu} because every open half-space through the origin contains either $\bm u_i$ or $-\bm u_i$ for some $i$. So the convergence results of Section \ref{sec:application-to-our-problem} apply to them, nevertheless, we solve it exactly here. The problem \eqref{eq:opt-problem} reduces to the following polar version
\begin{equation}\label{eq:finite-direction-optimization}
\inf_{Q\in \mathcal K^n_{\circ}} ~~ \sum_{i=1}^n \left[ a_i h_Q(\bm u_i)^\alpha + b_i h_Q(-\bm u_i)^\alpha \right] \quad \text{ s.t. } \quad \Phi(Q)\le 1.
\end{equation}

Let $U\in \R^{n\times n}$ be the matrix whose $i^{\text{th}}$ row is $\bm u_i^\top$, so that $U\bm x = (\inner{\bm u_1}{\bm x},\cdots, \inner{\bm u_n}{\bm x})$ as a column vector.  We show that the optimal body is the shifted zonotope
\begin{equation}\label{eq:solution-finite-support}
Q^* = U^{-1}\bm m^* + \sum_{i=1}^n r_i^*[-\bm v_i,\bm v_i]
\end{equation}
where $\bm v_i = U^{-1}\bm e_i$ and $\bm m^*, \bm r^*$ are appropriate shift and length parameters that make $\vol((Q^*)^\circ) = 1$. These parameters are as follows. Define
\[A_i\sett a_i^{1/(\alpha+1)},\quad B_i\sett b_i^{1/(\alpha+1)}, \quad S_i\sett A_i+B_i, \quad D_U\sett \frac{n!}{|\det U|}, \quad \kappa \sett D_U^{1/n} \prod_{j=1}^n S_j^{-(\alpha+1)/(\alpha n)}. \] Then \[m_i^* = \frac{(A_iS_i^{1/\alpha})^{-1} - (B_iS_i^{1/\alpha})^{-1}}{2\kappa}, \qquad r_i^* = \frac{(A_iS_i^{1/\alpha})^{-1} + (B_iS_i^{1/\alpha})^{-1}}{2\kappa}.\]

\cite[Proposition 3.1]{luo2020polar} identify the facet normals of the optimal body for an atomic measure. Here we extend such a result by using our KKT conditions to determine the translation and edge lengths of the optimal body as well.


\begin{thm}
\eqref{eq:solution-finite-support} gives the unique minimizer of \eqref{eq:finite-direction-optimization}. The corresponding primal optimizer is
\[
K^*=(Q^*)^\circ = \kappa \cdot \operatorname{conv} \set{ S_i^{1/\alpha} A_i\bm u_i, - S_i^{1/\alpha}B_i\bm u_i : 1\le i\le n}.
\]
\end{thm}

\begin{proof}
Let $E(h)$ be the objective above and denote $P=\set{\pm \bm u_i}_{i=1}^n\subseteq \Sn$. Applying the atomic derivative formula in Example \ref{ex:h^alpha} to the $2n$ directions $\set{\pm\bm u_i}_{i=1}^n$ with weights $a_i,b_i$ we get, \begin{equation}\label{eq:derivative-energy}
DE(h) = \alpha \sum_{i=1}^n\left[a_ih(\bm u_i)^{\alpha-1}\delta_{\bm u_i} + b_ih(-\bm u_i)^{\alpha-1}\delta_{-\bm u_i}\right].    
\end{equation} 
By a similar reasoning, \begin{equation}\label{eq:derivative-vol}
D\Phi(h) = -h^{-n-1}\sigma
\end{equation} where $\sigma$ is the unnormalized surface measure on $\mathbb S^{n-1}$. Note that $\inner{DE(h)}{h} = \alpha E(h)$ and $\inner{D\Phi(h)}{h} = -n\Phi(h)$. Now take the Lagrangian,
\[\mathcal L(h,\lambda,\nu) = E(h) + \lambda(\Phi(h)-1) - \inner{\nu}{h}.\]

We shall invoke the KKT conditions from Proposition \ref{prop:kktconditions} applied to the vector space of continuous functions $C(\Sn)$ with the open convex domain $\Omega = \set{f\in C(\Sn) : \min f >0}, \cF_0=E+\iota_{\mathscr F^n}, \cF_1=\Phi-1$. Note that Slater's condition holds for the constant support function $h_0\equiv R$ (on $\Sn$) where $R > \vol(\mathcal B_2^n)^{1/n}$ whence $E(h_0) < \infty$ and $\Phi(h_0) = \vol(R^{-1}\cB_2^n) = R^{-n}\vol(\mathcal B_2^n) < 1$. We try to find the minima $(h,\lambda,\nu)$ by solving the KKT conditions:
\begin{align}
h\in \mathscr F^n, h>0, \Phi(h)\le 1,\label{eq:feasible-atomic}\\
\lambda \ge 0, \nu \in (\mathscr F^n)^*, \label{eq:dualfeasible-atomic}\\
DE(h) + \lambda D\Phi(h) - \nu = 0,\label{eq:stationarity1-atomic}\\
\lambda(\Phi(h)-1) = 0,\label{eq:stationarity2-atomic}\\
\inner{\nu}{h} = 0.\label{eq:stationarity3-atomic}
\end{align}
\Cref{eq:stationarity1-atomic,eq:stationarity3-atomic} give $\inner{DE(h)}{h} = -\lambda \inner{D\Phi(h)}{h}$. Thus $\alpha E(h) = n\lambda \Phi(h)$. Since $h>0$, we have $E(h)>0$ and thus $\lambda>0$. By \eqref{eq:stationarity2-atomic}, $\Phi(h) = 1$. Recall that $\bm u_1,\cdots,\bm u_n$ form a basis. For every vector $\bm u\in \R^n$ let $c_i(\bm u) = \inner{\bm v_i}{\bm u} \in \R$ denote its coefficients in the expansion of this basis:
\begin{equation}\label{eq:basis-functions}\bm u = \sum_{i=1}^nc_i(\bm u)\bm u_i.\end{equation}
By sublinearity of $h$,
\begin{equation}\label{eq:h-compare}
h(\bm u) \le \sum_{i=1}^n h(c_i(\bm u) \bm u_i) = \sum_{i=1}^n \left[ c_i^+(\bm u) h(\bm u_i) + c_i^-(\bm u)h(-\bm u_i) \right]
\end{equation}
where $c_i^{\pm} = \max\{0,\pm c_i\}$. Define \[\tilde h(\bm u) \sett \sum_{i=1}^n \left[ c_i^+(\bm u) h(\bm u_i) + c_i^-(\bm u)h(-\bm u_i) \right].\]
Note that $c_i^+(\bm u) = \max\set{0,\inner{\bm v_i}{\bm u}} = h_{[0,\bm v_i]}(\bm u)$ is the support of the line segment between $0$ and $\bm v_i$ where $\bm v_i \sett U^{-1}\bm e_i$ form the dual basis. Similarly $c_i^- = h_{[-\bm v_i,0]}$. Therefore, $\tilde h \in \mathcal F$. Additionally note that $h,\tilde h$ match at the evaluation points:
\[\tilde h(\bm u_i) = h(\bm u_i), \qquad \tilde h(-\bm u_i) = h(-\bm u_i).\]
So we have \begin{equation}\label{eq:positive}\tilde h\ge h.\end{equation} This means that $\inner{DE(h)}{\tilde h} = \inner{DE(h)}{h}$ since $DE(h)$ is supported only in the directions $\set{\pm \bm u_i}_{i=1}^n$.  Dual feasibility gives \begin{equation}\label{eq:tightness}
0 \stackrel{\eqref{eq:dualfeasible-atomic}}{\le} \inner{\nu}{\tilde h} \stackrel{\eqref{eq:stationarity3-atomic}}{=} \inner{\nu}{\tilde h-h} = \lambda \inner{D\Phi(h)}{\tilde h-h} = \lambda \int h^{-n-1}\left(h-\tilde h\right)d\sigma \stackrel{\eqref{eq:dualfeasible-atomic}, \eqref{eq:positive}}{\le} 0.
\end{equation}
Since the integrand in the above expression is non-positive and continuous, due to \eqref{eq:positive}, we conclude that $h=\tilde h$ pointwise on the sphere. 

Translating the support function to its corresponding convex body, we get \[Q = \sum_{i=1}^n [-h(-\bm u_i)\bm v_i, h(\bm u_i)\bm v_i] = U^{-1} \prod_{i=1}^n [-h(-\bm u_i),h(\bm u_i)].\] 
This reduces to a finite dimensional problem over the variables $p_i=h(\bm u_i), q_i=h(-\bm u_i)$. Optimizing over these variables via KKT conditions subject to the polar-volume constraint gives the aforementioned solutions $\bm p^*=\bm r^*+\bm m^*, \bm q^*=\bm r^*-\bm m^*$. We shall briefly sketch the argument that verify this. Write $h=\sum_i(p_i c_i^++q_i c_i^-)$ and $s_i=p_i^{-1}+q_i^{-1}$. Dilating $h$ to $\Phi(h)^{1/n}h$ normalizes the polar volume, so the optimal $\bm p,\bm q$ minimize $\hat E(h)=\Phi(h)^{\alpha/n}E(h)$ over positive $(p_i,q_i)$. Using $\Phi(h)=|\det U|\prod_i s_i/n!$, stationarity in these variables gives
\[
a_i p_i^{\alpha+1}s_i=b_i q_i^{\alpha+1}s_i=\frac{E(h)}{n}.
\]
Thus $A_ip_i=B_iq_i$, and $s_i\propto S_i^{(\alpha+1)/\alpha}$. The constraint $\prod_i s_i=D_U$ fixes the proportionality constant to be $\kappa$, giving the stated parameters. These parameters result in the polar body in \eqref{eq:solution-finite-support} and the primal body is its polar. At the normalized minimizer, $D\hat E(h)=DE(h)+\frac{\alpha }{n}E(h)D\Phi(h)=\nu$. Since varying $p_i$ or $q_i$ varies $h$ in the direction $c_i^+$ or $c_i^-$ respectively, the same stationarity gives $\inner{\nu}{c_i^\pm} = 0$. We call this minimizer $h^*$.

Thus our proposed KKT certificate is $(h^*,\lambda^*,\nu^*) = (h^*, \alpha E(h^*)/n, DE(h^*)+ \alpha E(h^*) D\Phi(h^*)/n)$. It is straightforward to check that $p_i^*, q_i^* > 0$ so that $h|_{\Sn}>0$. All conditions \eqref{eq:feasible-atomic}--\eqref{eq:stationarity3-atomic} check out. Indeed, \eqref{eq:feasible-atomic}, \eqref{eq:stationarity1-atomic} hold simply by construction, \eqref{eq:stationarity2-atomic} holds because $\Phi(h^*)=1$, \eqref{eq:stationarity3-atomic} holds because $\inner{\nu}{c_i^\pm} = 0$ and $h^* = \sum(p_i^*c_i^++q_i^*c_i^-)$. Finally, the dual feasibility $\nu^*\in (\mathscr F^n)^*$ from \eqref{eq:dualfeasible-atomic} holds by almost the same argument used between Equations \eqref{eq:basis-functions} and \eqref{eq:positive} above. If $g\in \mathscr F^n$ is any support function, writing $\tilde g \sett \sum_i(c_i^+ g(\bm u_i) + c_i^- g(-\bm u_i))$ gives $\inner{\nu^*}{\tilde g} = 0$ since $\inner{\nu^*}{c_i^\pm} = 0$. But $\tilde g|_P = g|_P$ and $g\le \tilde g$ pointwise. Since $\operatorname{supp}DE(h^*) \subseteq P$, we get $\inner{DE(h^*)}{g-\tilde g} = 0$. Writing $\nu^*=DE(h^*) + \lambda^* D\Phi(h^*)$ gives $\inner{\nu^*}{g} = \inner{\nu^*}{g-\tilde g} = \lambda^* \inner{D\Phi(h^*)}{g-\tilde g} \stackrel{\eqref{eq:derivative-vol}}= \lambda^* \inner{(h^*)^{-n-1}\sigma}{\tilde g - g} \ge 0$, where the first equality used the fact that $\inner{\nu^*}{\tilde g}=0$, second equality used that $DE(h^*)$ is supported on $P$ and $g-\tilde g$ vanishes on $P$, and the last inequality on non-negativity used the fact that $g\le \tilde g$ pointwise and $(h^*)^{-n-1}\sigma$ is a strictly positive measure. This shows $\inner{\nu^*}{g}\ge 0$ for any $g\in \mathscr F^n$ whence \eqref{eq:dualfeasible-atomic} also holds. 

Conclude by Proposition \ref{prop:kktconditions}. $\lambda^*>0$ makes $E+\lambda^*\Phi$ strictly convex, giving uniqueness.
\end{proof}

Intuitively, this proof says that the evaluation directions need to be forced, and these are the only directions that are forced. Once the stretch along these directions are determined, all the other points in the boundary are simply an interpolation of these evaluation directions.

\begin{cor}
Suppose $\bm u_i=\bm e_i$ for $i=1,\cdots,n$. Then the polar and primal optimizers are
\[Q^* =  \bm m^*+\sum_{i=1}^n r_i^*[-\bm e_i,\bm e_i], \qquad K^* = \operatorname{conv} \set{ \frac{\bm e_i}{p_i^*}, -\frac{\bm e_i}{q_i^*} : i\in [n] }\]
where $p_i^* = r_i^* + m_i^*, q_i^* = r_i^* - m_i^*$.
\end{cor}

\begin{cor}
Suppose $\bm u_i=\bm e_i$ (standard basis) for $i=1,\cdots,n$. If $a_i=b_i\fa i\in[n]$ then $Q^*$ is a zero-centered zonotope. If, in addition, $a_i=b_i=(2n)^{-1}~\forall~i\in[n]$, then 
\[
Q^*=\vol(\cB_{\ell_1}^n)^{1/n}\mathcal B_{\ell_\infty}^n, \qquad K^*=\vol(\cB_{\ell_1}^n)^{-1/n}\mathcal B_{\ell_1}^n.
\]
\end{cor}

%% file: sec_Conclusion.tex
\section{Conclusion} \label{sec:conclusion}

We have developed a general framework for optimization over the space of convex bodies, with the main motivation being the problem of learning an optimal convex regularizer. In particular, we developed both optimality theory and analyzed the gradient flow dynamics of a generic variational problem when the objective and constraint set are characterized by functionals on the space of convex bodies. We applied this theory to the problem of learning an optimal convex regularizer. The two central steps were to reduce from $\mathscr K^n$ to $\body n$ for the convex regularizer selection problem, and then pass to the polar to reveal hidden convexity. This allowed us to exploit the aforementioned theory built for gradient flows. The KKT conditions were used in examples to demonstrate their power in characterizing the optimal body, with the $\ell_1$-ball being optimal for data characterized by an $\ell_p$-quasinorm being a notable example.

Below we discuss important future directions:
\begin{itemize}
    \item \textbf{Design principles for the optimal convex body:} The two examples studied in Section \ref{sec:optimal-convex-reg-examples} illustrate the information about the optimizer hidden in the dual certificates. For example, they reveal the extreme points of the optimal convex body studied in the problem. In particular, the support and symmetry of the initial data distribution hold information about the geometric and combinatorial structure of the optimal regularizer. A natural next step is to explore along this direction of making suitable symmetry assumptions to recover specific convex regularizers. It is important to understand how certain symmetries and geometric properties of the data distribution characterize the optimal convex regularizer.
    \item \textbf{Distributional shifts:} Though our example concerning $\mu(\bm x) \propto\exp(-\|\bm x\|_p)$ for $p\in(0,1)$ results in polyhedral optimizers and \cite[Section 4.5]{leong2025optimal} implies that the optimal convex body $\vol(K)^{-1/n}K$ is optimal for the measure $\mu(\bm x) \propto \exp(-\|\bm x\|_K)$ if $K$ is convex, it is important to ask about quantitative rates about the change in the optimal body as the $K$ in question becomes smoother or transitions from a star body to a smooth convex body.
    \item \textbf{Statistical learning aspects:} A related statistical question arises when a finite sample replaces the population distribution. Given a target accuracy in Hausdorff distance, how many samples suffice to estimate the optimal convex body with high probability? Also, what lower bounds on the sample size are unavoidable?
    \item \textbf{Geometry of bodies along the flow and transport:} Finally, the convergence results raise questions about the geometry of the bodies \textit{along} the flow. For example, if the initial body is a polytope and the limiting optimizer has a smooth boundary, under what conditions does the evolving body acquire a smooth boundary in finite time? Optimal transport offers another way to study this evolution. For example, associating each body with its normalized volume measure would allow us to compare the path induced by the gradient flow with an optimal transport interpolation between the initial and optimal bodies. It would then be natural to compare the metrics $\delta_p$ against these transport metrics, up to appropriate normalization.
\end{itemize}

\section*{Acknowledgements}
We acknowledge the usage of ChatGPT and Claude for help with literature search and editing. Part of this research was conducted while NM was visiting the National University of Singapore (NUS).  NM gratefully acknowledges financial support from Jonathan Scarlett and Yong Sheng Soh for visiting NUS.  YS acknowledges support from Ministry of Education Tier 1 Grants 24-0752-A0001 and 26-0787-A0001.